\documentclass[12pt,reqno]{amsart}
\usepackage{amsmath}
\usepackage{amssymb}
\usepackage{amsthm}
\usepackage{mathrsfs}
\usepackage{latexsym}
\usepackage{xcolor}
\usepackage{comment}

\newtheorem{theorem}{Theorem}[section]
\newtheorem{lemma}{Lemma}[section]
\newtheorem{corollary}{Corollary}[section]
\newtheorem{remark}{Remark}[section]
\newtheorem{definition}{Definition}[section]
\newtheorem{proposition}{Proposition}[section]
\newtheorem{example}{Example}[section]
\newtheorem{assumption}{Assumption}[section]
\numberwithin{equation}{section}
\newcommand{\bth}{\begin{theorem}}
\newcommand{\ethe}{\end{theorem}}

\newcommand{\bre}{\begin{remark}}
\newcommand{\ere}{\end{remark}}

\newcommand{\ble}{\begin{lemma}}
\newcommand{\ele}{\end{lemma}}
\newcommand{\bde}{\begin{definition}}
\newcommand{\ede}{\end{definition}}

\newcommand{\bco}{\begin{corollary}}
\newcommand{\eco}{\end{corollary}}

\newcommand{\bpr}{\begin{proposition}}
\newcommand{\epr}{\end{proposition}}

\newcommand{\bexer}{\begin{exercise}}
\newcommand{\eexer}{\end{exercise}}
\newcommand{\breh}{\begin{hint}}
\newcommand{\ereh}{\end{hint}}

\newcommand{\halmos}{\hfill \qed}

\newcommand{\bexam}{\begin{example}}
\newcommand{\eexam}{\end{example}}

\newcommand{\pr} {{\bf Proof.}}

\newcommand{\bfi}{\begin{fig}}
\newcommand{\efi}{\end{fig}}

\newcommand{\beao}{\begin{eqnarray*}}
\newcommand{\eeao}{\end{eqnarray*}\noindent}

\newcommand{\beam}{\begin{eqnarray}}
\newcommand{\eeam}{\end{eqnarray}\noindent}
\newcommand{\E}{\mathbf{E}}
\newcommand{\PP}{\mathbf{P}}

\newcommand{\xto}{x\to\infty}

\newcommand{\bF}{\overline{F}}

\newcommand{\bG}{\overline{G}}
\newcommand{\bV}{\overline{V}}

\newcommand{\bbr}{{\mathbb R}}

\newcommand{\bbn}{{\mathbb N}}

\newcommand{\vep}{\varepsilon}

\allowdisplaybreaks[1]

\begin{document}
\title[Asymptotic estimates for a nonstandard bidimensional risk model]{Asymptotic estimations for a nonstandard bidimensional risk model with dependent claims and constant interest force}

\author[D.G. Konstantinides, C.D. Passalidis]{Dimitrios G. Konstantinides, Charalampos  D. Passalidis} 

\address{Dept. of Statistics and Actuarial-Financial Mathematics,
University of the Aegean, Karlovassi, GR-83 200 Samos, Greece}
\email{konstant@aegean.gr}
\email{sasd24009@sas.aegean.gr.}

\date{{\small \today}}

\begin{abstract}
In this paper we extend the results from Yang and Li(2017), in a non-renewal driven bivariate risk model. Concretely, we are interested in the asymptotic behavior of the joint tail of  the discounted aggregate claims over finite and infinite time horizon in a bivariate risk model with two arbitrarily dependent counting processes. We additionally suppose that the 
sequences of claim-sizes of the two lines of business are independent, but each of them contains weakly dependent terms. In our main results, on finite and on infinite time horizon, we assume two (different) general conditions for the counting processes, that are satisfied by a wide spectrum of processes beyond the renewal ones. In the case of finite horizon we suppose 
that the claim distributions from the two lines, belong to the subexponential distribution class, while in the case of infinite horizon we restrict ourselves to the 
consistently varying and positively decreasing  distribution class. More explicit asymptotic expressions are derived in the case we are restricted on the regularly varying class for the
claim distributions. We note that our results indicate the presence of multivariate non-linear single big jump principle for the discounted aggregate claims.
\end{abstract}

\maketitle
\vspace{3mm}




\section{Introduction} \label{sec.CKP29V.1}

\subsection{Model description} \label{subsec.CKP29V.1.1}

We consider an insurer, who operates two lines of business. We suppose that the claim arrivals into these two lines happen at moments $\{\tau_i\,,\;i \in \bbn\}$ and $\{\sigma_j\,,\;j \in \bbn\}$, respectively, with $\tau_0= \sigma_0=0$. The sequences of the arrival times are non-degenerate to zero, and constitute the counting processes
\beam \label{eq.CKP29V.1.1}
N(t):= \sup \{n \in \bbn\;:\;\tau_n \leq t\}\,, \quad M(t):= \sup \{n \in \bbn\;:\;\sigma_n \leq t\}\,,
\eeam
with finite expectations
\beao
\nu(t) := \E[N(t)] = \sum_{i=1}^{\infty} \PP[\tau_i \leq t]\,, \quad \mu(t) := \E[M(t)] = \sum_{j=1}^{\infty} \PP[\sigma_j \leq t]\,,
\eeao
for any $t\geq 0$. Their joint expectation is given by
\beao
\E[N(s)\,M(t)] := \sum_{i=1}^{\infty} \sum_{j=1}^{\infty} \PP(\tau_i \leq s\,,\;\sigma_j \leq t)\,,
\eeao
for $s,\,t\geq 0$. The claims of the two lines of business are described by two sequences of non-negative random variables (r.v.s) $\{X_i,\,i \in \bbn\}$, $\{Y_j,\,j \in \bbn\}$, and each of them has identically distributed (but not necessarily independent) terms, 
with common distributions denoted by $F$ and $G$, respectively. We suppose that the insurer invests the surplus of each line, and obtains constant interest forces 
$r_1\,,\;r_2 \geq 0$, respectively. Hence the discounted aggregate claims up to moment $t \geq 0$, are described by the relation
\beam \label{eq.CKP29V.1.2}
{\bf D}(t) = \begin{pmatrix}   
D_1(t)  \\ 
D_2(t) 
\end{pmatrix}= \begin{pmatrix}   
\sum_{i=1}^{N(t)} X_i\,e^{-r_1\,\tau_i}  \\[2mm]
\sum_{j=1}^{M(t)} Y_j\,e^{-r_2\,\sigma_j}   
\end{pmatrix}\,, 
\eeam
where in case $t=\infty$, we write $N(\infty)= M(\infty)= \infty$.

Our main goal in this paper is to provide precise asymptotic expressions for the joint tail of the discounted aggregate claims
\beam \label{eq.CKP29V.1.3}
\PP[D_1(t)  > x\,,\;D_2(t) >y]\,,
\eeam
as $(x,\,y) \to (\infty,\,\infty)$, for both finite horizon $t=T$, and infinite  horizon $t= \infty$ as well, under some heavy tail conditions for $F$ and $G$. Our 
results indicate the presence of the multivariate non-linear single big jump principle for \eqref{eq.CKP29V.1.3}, see \cite{konstantinides:passalidis:2025c} for more 
details about this principle.

For the main results of the paper, we make use of the following assumption, for which we need the following sequences of non-negative and non-degenerate to zero r.v.s $\{\Theta_i := \tau_i - \tau_{i-1}\,,\; i \in \bbn\}$ and  $\{\Delta_j := \sigma_j - \sigma_{j-1}\,,\; j \in \bbn\}$. These r.v.s represent the interarrival times between two successive claim arrivals in each line of business.

\begin{assumption} \label{ass.CKP29V.1.1}
The sequences $\{X_i,\,i \in \bbn\}$, $\{Y_j,\,j \in \bbn\}$, $\{\Theta_i,\,\Delta_j\,,\; i,\,j \in \bbn\}$ are mutually independent.
\end{assumption}

\subsection{Brief review} \label{subsec.CKP29V.1.2}

The bivariate or more general the multivariate risk models represent a more realistic environment for the actuarial applications, since the modern insurance companies 
operate multiple portfolios, in order to be competative. Hence, during the recent years, many researchers focused their interest in the study of multvariate risk models, mainly in cases of claims with heavy tails.

Due to the good properties of the multivariate regularly varying distributions, and the multivariate subexponentiality, introduced by \cite{samorodnitsky:sun:2016}, 
most of the papers studied various risk models, with main characteristic the asymptotic estimations for the discounted aggregate claims (or, the, derived by them, 
ruin probabilities) under the multivariate linear single big jump principle, see for example in \cite{konstantinides:li:2016}, \cite{konstantinides:passalidis:2025j}, 
\cite{passalidis:2025}, \cite{konstantinides:passalidis:xu:2026} and \cite{yuan:lu:fu:2025}, among others. However, the linear approximation of the single big jump principle has 
the drawback that the convergences of initial capitals of the lines of business are made with respect to $l_i\,x \to \infty$, as $x \to \infty$, with some $l_i \in (0,\,\infty)$, and therefore in many cases is underestimated the difference in the heaviness of the claim-tails among these lines. Additionally, in the case of multivariate subexponentiality the events joint tail are not possible to be studied, see in \cite[Rem. 2.2]{konstantinides:passalidis:2025g}. So, in the case of non-homogeneous risks, while the linear approach works well in the study of sum and maximum sets, since appears domination of the heaviest tail, it does not happen the same for the joint tail. The case of standard multivariate regular variation, that can study also the case of joint tail, often also fails in the case of non-homogeneous risks, with the compensation by hidden multivariate regular variation, see \cite{resnick:2024}, to seem very restrictive for the applications of the risk theory.

Finally, we note that in most of the previous papers there is common counting process, something that can be restrictive in actuarial practice, especially for the cases we use simple common processes, like the renewal driven risk models.

All these facts, bring into the light the need of additional study of risk models, whose claims satisfy, in some sense, the multivariate non-linear single big jump principle, see in \cite{konstantinides:passalidis:2025c}, \cite{konstantinides:passalidis:2025g} for more details about the distinction of these principles.

The study of various risk models with emphases in the joint tail, and furthermore in the multivariate non-linear single big jump, is often restricted to bidimensional set up, because of the more difficulty in the non-linear in comparison with the linear approach. In \cite{chen:yuen:ng:2011}, \cite{chen:wang:wang:2013}, \cite{gao:yang:2014}, \cite{jiang:wang:chen:xu:2015}, \cite{liu:geng:man:liu:2023} we have only some of the contributions to this direction. However, the previous papers focus mostly in the cases when either there is common counting process for the two lines of business, or there is independence between these two counting processes.

In the pioneer paper \cite{yang:li:2017} we find for first time the study of risk model  \eqref{eq.CKP29V.1.2}, with $r_1=r_2$, under the Assumption \ref{ass.CKP29V.1.1}, considering that the $\{(\Theta_i,\,\Delta_i)\,,\;i \in \bbn\}$ are independent and identically distributed (i.i.d.) copies of $(\Theta,\,\Delta)$, however each pair has arbitrarily dependent components. Hence, in their model the processes $\{N(t),\,t \geq 0 \},\,\{M(t),\,t \geq 0\}$ are two arbitrarily dependent (homogeneous) renewal processes. Under some heavy-tailed conditions on  distributions $F,\,G$ were provided asymptotic expressions over finite and over infinite horizon for probability \eqref{eq.CKP29V.1.3}, as also for a concrete kind of ruin probability, when the $\{X_i,\,i \in \bbn\}$, $\{Y_j,\,j \in \bbn\}$, are both i.i.d. sequences. 

The goal of this paper is to provide the same asymptotic expressions over finite and over infinite horizon, under much more general conditions, as follows:
\begin{enumerate}
\item
The counting processes  $\{N(t),\,t \geq 0 \},\,\{M(t),\,t \geq 0\}$ are not necessarily renewal, see Assumptions \ref{ass.CKP29V.3.1} and \ref{ass.CKP29V.3.2} below.
\item
The sequences $\{X_i,\,i \in \bbn\}$, $\{Y_j,\,j \in \bbn\}$, contain weakly dependent terms,  see Assumptions \ref{ass.CKP29V.2.1} and \ref{ass.CKP29V.2.2}. 
\item
In case $t=\infty$, we widen the spectrum, of the distribution classes of $F$ and $G$, namely from extended regularly varying claims we move to $\mathcal{C}\cap \mathcal{P_D}$.
\end{enumerate}

We also note that several papers follow \cite{yang:li:2017} in the study of model \eqref{eq.CKP29V.1.2} and its generalizations. In \cite{chen:yang:jiang:2019} was 
studied the ruin probability in the renewal risk model  \eqref{eq.CKP29V.1.2}, and derived local uniform, with respect to time, asymptotic estimations when the 
$\{X_i,\,i \in \bbn\}$, $\{Y_j,\,j \in \bbn\}$ satisfy a general dependence structure, see Assumption \ref{ass.CKP29V.2.2}, while the $F$ and $G$ belong to class $\mathcal{D} \cap \mathcal{L}$. In \cite{chen:li:cheng:2023} was examined another kind of ruin probability in a more general model than that of \eqref{eq.CKP29V.1.2}, where the logarithmic returns of the insurer's investment portfolios are described by two non-negative L\'{e}vy processes. Further, in that paper, was assumed that the $\{X_i,\,i \in \bbn\}$, $\{Y_j,\,j \in \bbn\}$ are weakly dependent, see Assumption \ref{ass.CKP29V.2.1}, and the processes $\{N(t),\,t \geq 0 \},\,\{M(t),\,t \geq 0\}$ satisfy some convergence conditions and an assumption slightly more strict than Assumption \ref{ass.CKP29V.3.1} below. In \cite{cheng:yang:wang:2020} was considered a variation of  Assumption \ref{ass.CKP29V.1.1}, where the $\{(X_i,\,Y_i)\,,\; i \in \bbn\}$ are i.i.d. copies of $(X,\,Y)$, that contains some weak form of dependence, where was used also Assumption \ref{ass.CKP29V.3.1} for the study of three kinds of ruin probabilities over finite horizon.   
 
The rest of this paper is organized as follows. In Section 2, we introduce the necessary preliminary concepts about the heavy-tailed distributions and the dependence 
structures we shall need later. In Section 3, we give two main results for asymptotic estimation of probability \eqref{eq.CKP29V.1.3} over finite and over infinite horizon, having already introduced the necessary assumption on the counting processes. In Section 4, we provide the proofs of the main results after some preliminary lemmas.

\section{Preliminaries} \label{sec.CKP29V.2}

Hereafter, all the limit relations hold as $(x,\,y) \to (\infty,\,\infty)$, except is written differently. For two positive univariate functions $\widetilde{f},\,\widetilde{g}$, we write $\widetilde{f}(x)=O[\widetilde{g}(x)]$, as $\xto$, if it holds
\beao
\limsup_{\xto} \dfrac{\widetilde{f}(x)}{\widetilde{g}(x)}< \infty \,,
\eeao
while we write $\widetilde{f}(x) \asymp \widetilde{g}(x)$, as $\xto$, if hold both $\widetilde{f}(x)=O[\widetilde{g}(x)]$ and $\widetilde{g}(x)=O[\widetilde{f}(x)]$, as $\xto$. Further, we write $\widetilde{f}(x) \sim \widetilde{g}(x)$, as $\xto$, if it holds
\beao
\lim_{\xto} \dfrac{\widetilde{f}(x)}{\widetilde{g}(x)} =1 \,.
\eeao
For two positive bivariate functions $f(\cdot,\,\cdot),\;g(\cdot,\,\cdot)$, we write $f(x,\,y) \sim g(x,\,y)$, $f(x,\,y) \lesssim g(x,\,y)$ and $f(x,\,y)=o[g(x,\,y)]$, if hold
\beao
\lim \dfrac{f(x,\,y)}{g(x,\,y)} = 1\,, \qquad \limsup \dfrac{f(x,\,y)}{g(x,\,y)} \lesssim 1\,, \qquad \limsup \dfrac{f(x,\,y)}{g(x,\,y)} = 0\,,
\eeao
respectively. We write $f(x,\,y)=O[g(x,\,y)]$, if it holds
\beao
\limsup \dfrac{f(x,\,y)}{g(x,\,y)} < \infty\,,
\eeao
and  $f(x,\,y)\asymp g(x,\,y)$, if hold simultaneously $f(x,\,y)=O[g(x,\,y)]$ and $g(x,\,y)=O[f(x,\,y)]$. Further, for two positive trivariate functions $f^*(\cdot,\,\cdot;\,\cdot)$, $g^*(\cdot,\,\cdot;\,\cdot)$, we say that it holds 
\beao
f^*(x,\,y;\,z) \sim g^*(x,\,y;\,z)\,,
\eeao 
uniformly for $z \in \Delta$, where $\Delta \neq \emptyset$, if it holds
\beao
\lim \sup_{z \in \Delta} \left| \dfrac{f^*(x,\,y;\,z)}{g^*(x,\,y;\,z)} -1\right| =0\,. 
\eeao
Similar notations of uniform convergence used also for bivariate positive functions.
For two real numbers $a,\,b$ we denote $a\wedge b :=\min \{a,\,b\}$ and $a\vee b :=\max\{a,\,b\}$. Finally  for a random variable $Z$ that follows distribution $V$, we write $Z \stackrel{d}{\sim} V$.

\subsection{Heavy-tailed distributions} \label{subsec.CKP29V.2.1}

In this subsection we provide the necessary preliminary definitions for the heavy tailed distribution classes, that we use later. We note that the main classes we use, 
are the classes $\mathcal{S}$ and $\mathcal{C} \cap \mathcal{P_D}$, for our results over finite and infinite horizon correspondingly, while the class of the regular 
variation gives some corollaries with more explicit asymptotic formulas. In all this section, the distribution $V$, is such that $\bV(x) >0$, for any $x \in \bbr$. For the sake of compactness, the following definitions will be presented for distributions with support in $\bbr_{+}:=[0,\,\infty)$, instead of $\bbr$.

At first, we say that a distribution $V$ has heavy tail, symbolically $V \in \mathcal{K}$, if for any $\vep >0$ it holds
\beao
\int_0^{\infty} e^{\vep\,y}\,V(dy) = \infty\,.
\eeao

We say that the distribution $V$ belongs to the class of long tailed distributions, symbolically $V \in \mathcal{L}$, if for any (or, equivalently, for some) $a>0$ it holds 
\beao
\lim_{\xto} \dfrac{\bV(x-a)}{\bV(x)}=1\,,
\eeao
while we say that the distribution $V$ belongs to the class of subexponential distributions, symbolically $V \in \mathcal{S}$, if for any (or, equivalently, for some) integer $n \geq 2$ it holds 
\beao
\lim_{\xto} \dfrac{\overline{V^{n*}}(x)}{\bV(x)}=n\,,
\eeao
where $V^{n*}$ represents the $n$-th fold convolution of $V$ with itself. These classes were introduce by \cite{chistyakov:1964}, where was also proved the inclusion $\mathcal{S} \subsetneq \mathcal{L} \subsetneq \mathcal{K}$. We refer the reader to \cite{athreya:ney:1972}, \cite{asmussen:albrecher:2010}, \cite{foss:korshunov:zachary:2013}, for properties and applications of these classes.

The class $\mathcal{D}$ of the dominatedly varying distributions, was introduced on \cite{feller:1969}. We say that $V \in \mathcal{D}$, if for any (or, equivalently, for some) $b\in (0,\,1)$ it holds
\beao
\limsup_{\xto} \dfrac{\bV(b\,x)}{\bV(x)} < \infty\,.
\eeao
It is well-known that $\mathcal{D} \subsetneq \mathcal{K}$, $\mathcal{D} \not\subseteq \mathcal{S}$, $\mathcal{S} \not\subseteq \mathcal{D}$, and 
$\mathcal{D} \cap \mathcal{S} \equiv \mathcal{D} \cap \mathcal{L}\neq\emptyset$, see \cite{goldie:1978}. We say that the distribution $V$ belongs to the class of consistently varying distributions, symbolically $V \in\mathcal{C}$, if
\beao
\lim_{b \uparrow 1} \limsup_{\xto} \dfrac{\bV(b\,x)}{\bV(x)} =1\,.
\eeao
The class of regularly varying distributions is defined as follows: We say that $V$ is regularly varying with index $\alpha \geq 0$, symbolically $V \in \mathcal{R}_{-\alpha}$, if for any $b > 0$ it holds
\beao
\lim_{\xto} \dfrac{\bV(b\,x)}{\bV(x)} = b^{-\alpha}\,.
\eeao
Often we write $\mathcal{R}_0$ instead of $\mathcal{R}_{-0}$, and this class is called class of slowly varying distributions (or, and more generally functions). For example a positive function $f$ with support on $(0,\,\infty)$, we write $f \in \mathcal{R}_0$, if for any $b>0$ it holds
\beao
\lim_{\xto} \dfrac{f(b\,x)}{f(x)} = 1\,.
\eeao

The following inclusions are well known, see for example  \cite[Ch. 2]{leipus:siaulys:konstantinides:2023} 
\beao
\mathcal{R}:= \mathcal{R}_0 \bigcup \left(\bigcup_{\alpha >0} \mathcal{R}_{-\alpha} \right) \subsetneq \mathcal{C} \subsetneq \mathcal{D} \cap \mathcal{L} \subsetneq \mathcal{S} \subsetneq \mathcal{L} \subsetneq \mathcal{K}\,.
\eeao

Class $\mathcal{P_D}$ of positively decreasing distributions was introduced on \cite{haan:resnick:1984} and is such that $\mathcal{P_D}\cap \mathcal{K} \neq \emptyset$, $\mathcal{P_D}\cap \mathcal{K}^c \neq \emptyset$. We say that $V \in \mathcal{P_D}$, if for any (or, equivalently, for some) $v>1$, it holds
\beao
\limsup_{\xto} \dfrac{\bV(v\,x)}{\bV(x)} < 1\,.
\eeao 
Class $\mathcal{P_D}$ in combination with the previous distribution classes plays crucial role in applied probability, especially in problems of infinite horizon, see for example  \cite{konstantinides:tang:tsitsiashvili:2002}, \cite{hao:tang:2008}.

In \cite{konstantinides:tang:tsitsiashvili:2002} was also introduced the class $\mathcal{A}:=\mathcal{S} \cap \mathcal{P_D}$, that contains all the usually used 
subexponential distributions, see \cite[Sec. 2]{tang:2006a}. We refer the reader to \cite{bardoutsos:konstantinides:2011}, \cite{konstantinides:passalidis:2024d} for 
more properties of $\mathcal{P_D}$ and related classes. Here, we shall use class $\mathcal{C} \cap \mathcal{P_D}$ for the result in the infinite horizon, and we note that class $\mathcal{C} \cap \mathcal{P_D}$ is slightly smaller than $\mathcal{C}$.

The following indexes indicate the symmetry between classes $\mathcal{D}$ and $\mathcal{P_D}$. The upper and lower Matuszewska indexes are correspondingly defined as follows:
\beam \label{eq.CKP29V.2.1}
J_V^+:=\lim_{v \to \infty} \dfrac{\log \bV_*(v)}{\log v}\,,\qquad J_V^-:=\lim_{v \to \infty} \dfrac{\log \bV^*(v)}{\log v}\,,
\eeam
where
\beao
\bV^*(v) = \limsup_{x \to \infty} \dfrac{\bV(v\,x)}{\bV(x)}\,,\qquad  \bV_*(v) = \liminf_{x \to \infty} \dfrac{\bV(v\,x)}{\bV(x)}\,. 
\eeao
From \eqref{eq.CKP29V.2.1} it is easy to see that for any distribution $V$, with $\bV(x) >0$ for any $x \in \bbr$, it holds $0 \leq J_V^- \leq J_V^+ \leq \infty$.

It is well-known that $V \in \mathcal{D}$ if and only if $J_V^+ <  \infty$, $V \in \mathcal{P_D}$ if and only if $J_V^- >0$,  and if $V \in \mathcal{R}_{-\alpha}$, with $\alpha \geq 0$, then $J_V^+ =J_V^- = \alpha$, see \cite[Subsec. 2.4]{leipus:siaulys:konstantinides:2023}. Under the previous characterization, we can see that the following inclusions are true:
\beam \label{eq.CKP29V.2.2}
\bigcup_{\alpha >0} \mathcal{R}_{-\alpha}  \subsetneq \mathcal{C} \cap \mathcal{P_D} \subsetneq \mathcal{D} \cap \mathcal{A} \subsetneq \mathcal{A} \subsetneq \mathcal{S} \subsetneq \mathcal{L} \subsetneq \mathcal{K}\,.
\eeam

We also notice that if $V \in \mathcal{D} \cap \mathcal{P_D}$ (that means $0 < J_V^- \leq J_V^+ < \infty$), then from \cite[Prop. 2.2.1]{bingham:goldie:teugels:1987} 
we obtain that for any pair $(p_1,\,p_2)$, that satisfy the inequalities $0< p_1 <J_V^- \leq J_V^+ < p_2 < \infty$, there exist constants $C_i>0$ and $D_i>0$, for $i =1,\,2$, such that it holds
\beam \label{eq.CKP29V.2.3}
\dfrac{\bV(y)}{\bV(x)} \geq C_1\,\left( \dfrac xy \right)^{p_1}\,,
\eeam
for any $x\geq y\geq D_1$, and also holds
\beam \label{eq.CKP29V.2.4}
\dfrac{\bV(y)}{\bV(x)} \leq C_2\,\left( \dfrac xy \right)^{p_2}\,,
\eeam
for any $x\geq y\geq D_2$. Furthermore, form \cite[Lem. 3.5]{tang:tsitsiashvili:2003}, for any $p_2 > J_V^+$ we find
\beam \label{eq.CKP29V.2.5}
\lim_{x \to \infty}\dfrac{x^{-p_2}}{\bV(x)} =0\,.
\eeam

\subsection{Dependence structures} \label{subsec.CKP29V.2.2}

In this section we provide the two dependence structures, that will be used later for the sequences of claim sizes $\{X_i,\,i \in \bbn\}$ and $\{Y_j,\,j \in \bbn\}$. Additionally, we assume for any $i \in \bbn$, that $Z_i \stackrel{d}{\sim} V_i$, with $\bV_i(x) > 0$, for any $x \in \bbr$ and the distributions $V_i$ have support on $\bbr_{+}$.

The following dependence structure was introduced by \cite{yang:liu:huang:ma:2012}, and represents generalization of the dependence from \cite{ko:tang:2008} over infinite number of summands. 

\begin{assumption} \label{ass.CKP29V.2.1}
For $\{Z_i,\,i \in \bbn\}$, let suppose that there exist constants $C>0$, $x_0 >0$, such that it holds
\beam \label{eq.CKP29V.2.6}
\sup_{x \geq x_0} \sup_{n\geq 2} \sup_{a \in [x_0,\,x]} \dfrac{\PP\left(\sum_{i=1}^{n-1} Z_i > x-a\;|\;Z_{n} = a \right)}{\PP\left(\sum_{i=1}^{n-1} Z_i > x-a \right)} \leq C\,.
\eeam
\end{assumption}

\bre \label{rem.CKP29V.2.1}
Assumption \ref{ass.CKP29V.2.1}, obviously contains the independence as special case. Further, it provides a general framework, that is consistent with the 
'insensitivity' of the single big jump principle, in the case of class $\mathcal{S}$ for the distributions of $\{Z_i,\,i \in \bbn\}$, see for example 
\cite[Th. 3.1]{ko:tang:2008} and relative discussions in \cite{jiang:gao:wang:2014}. For these reasons Assumption \ref{ass.CKP29V.2.1} and its extensions have been 
used in many applications of risk theory, see for example \cite{zhang:cheng:2017},  \cite{chen:li:cheng:2023}, \cite{chen:wang:cheng:yan:2023}. We refer the reader to \cite{yang:liu:huang:ma:2012} and \cite{zhang:cheng:2017}, for examples that satisfy Assumption \ref{ass.CKP29V.2.1}.
\ere

As we mentioned previously, in the case of infinite horizon, we restrict the distributions $F$ and $G$ into $\mathcal{C} \cap \mathcal{P_D}$ instead of $\mathcal{S}$. However, the dependence structure we use is more general than that of Assumption \ref{ass.CKP29V.2.1}. This dependence was introduces by \cite{geluk:tang:2009}, and was used in several papers in risk theory and in risk management, see for example \cite{chen:wang:wang:2013}, \cite{li:2013}, \cite{chen:liu:2022}, among many others. 

\begin{assumption} \label{ass.CKP29V.2.2}
Let suppose that the $\{Z_i,\,i \in \bbn\}$ are pairwise tail asymptotically independent, namely for any $i,\,j \in \bbn$, with  $i \neq j$, it holds
\beam \label{eq.CKP29V.2.7}
\lim_{x_i \wedge x_j \to \infty} \PP\left( Z_i > x_i\;|\;Z_j > x_j \right)=0\,.
\eeam
\end{assumption}

\section{Main results} \label{sec.CKP29V.3}

In this section we present our main results on the asymptotic estimations of \eqref{eq.CKP29V.1.3} over finite and infinite horizon. For the case of finite horizon we 
make use of the following general assumption for the counting processes.

\begin{assumption} \label{ass.CKP29V.3.1}
The pair $(\Theta_1,\,\Delta_1)$ is independent of $\{(\Theta_i,\,\Delta_j)\,,\; i,\,j \in \bbn\,,\; i \wedge j \geq 2\}$.
\end{assumption}

\bre \label{rem.CKP29V.3.1}
The dependence between insurer's two lines of business, under the Assumption \ref{ass.CKP29V.1.1}, stems only from the dependence between the $\{N(t),\,t \geq 0 \}$  
and $\{M(t),\,t \geq 0\}$. We can see that Assumption \ref{ass.CKP29V.3.1} contains a wide spectrum of counting processes, that escape from the common renewal driven 
models and their generalizations. A such kind of assumption, in combination with its application on risk models with subexponential claims, was used first time by \cite{cheng:yang:wang:2020}, while another one slightly stricter (with additional condition that the $\Theta_1,\,\Delta_1$ are independent each other), was used by \cite{chen:li:cheng:2023}. The assumption  that the $(\Theta_1,\,\Delta_1)$ are independent of the rest of inter-arrival times is very helpful from mathematical point of view, since it permits application of the Kesten inequalities. From practical point of view Assumption \ref{ass.CKP29V.3.1} is not much stricter than the condition: ' $\{\Theta_i,\,\Delta_j\,,\; i,\,j \in \bbn\}$ are arbitrarily dependent', since in many practical cases at the beginning of the operation of a portfolio-company, appear some 'irregularities' that get smoothed out by the time, see for example \cite{wang:2008}, for relative discussions. 
\ere

\bre \label{rem.CKP29V.3.2}
For the conditions of our first result, we need the following delayed counting processes:
\beam \label{eq.CKP29V.3.1}
N^*(t):=\sup \left\{n \in \bbn\;:\; \sum_{i=1}^n \Theta_{i+1} \leq t \right\}\,,\;M^*(t):=\sup \left\{n \in \bbn\;:\; \sum_{j=1}^n \Delta_{j+1} \leq t \right\}\,, 
\eeam
for any $t\geq 0$.
\ere

The following theorem is the first main result of this paper.

\bth \label{th.CKP29V.3.1}
Let consider the discounted aggregate claims of relation \eqref{eq.CKP29V.1.2}. We suppose that Assumptions \ref{ass.CKP29V.1.1} and \ref{ass.CKP29V.3.1} are valid. For each of the sequences of r.v.s $\{X_i,\,i \in \bbn\}$ and $\{Y_j,\,j \in \bbn\}$ Assumption \ref{ass.CKP29V.2.1} is satisfied. If $F,\,G \in \mathcal{S}$, then for any fixed $T>0$, such that $\E[N(T)\,M(T)] >0$, and for which there exists some constant $\beta = \beta(T) >0$, such that the conditions
\beam \label{eq.CKP29V.3.2}
\E \left[e^{\beta\,N^*(T)}\right] < \infty\,,\;\E \left[e^{\beta\,M^*(T)}\right] < \infty \,,
\eeam
are true, we obtain 
\beam \label{eq.CKP29V.3.3}
\PP[D_1(T)  > x\,,\;D_2(T) >y] \sim \int_0^T \int_0^T \bF(x\,e^{r_1\,s})\,\bG(y\,e^{r_2\,t})\,\E[N(ds)\,M(dt)]\,.
\eeam   
\ethe

\bre \label{rem.CKP29V.3.3}
Theorem \ref{th.CKP29V.3.1} generalizes \cite[Lem. 4.4]{yang:li:2017}, with respect to (1) and (2) of subsection \ref{subsec.CKP29V.1.2}. The generalization with respect to (1), follows from the fact that if the $\{N(t),\,t \geq 0 \}$ and $\{M(t),\,t \geq 0\}$ are renewal processes, then \eqref{eq.CKP29V.3.2} holds, and even for any fixed $T >0$. Indeed, in such a case, the delayed processes from \eqref{eq.CKP29V.3.1} should be also renewal, and hence from \cite{stein:1946} they have analytic moment-generating functions near to zero.
\ere

\bre \label{rem.CKP29V.3.4}
The statement of Theorem \ref{th.CKP29V.3.1} depends on the choice of $T>0$, due to condition \eqref{eq.CKP29V.3.2}. However, as comes from the Remark \ref{rem.CKP29V.3.3}, in the renewal case, through Stein's lemma, for any fixed $T > 0$, there exists a properly chosen $\beta=\beta(T) >0$, such that relation \eqref{eq.CKP29V.3.2} holds. 

The same is true for some quasi-renewal processes (i.e counting processes with identical but not necessarily independent inter-arrival times), see \cite[Lem. 2.2]{wang:cheng:2011}, and for some non-homogeneous renewal processes  (i.e. independent but not necessarily identical inter-arrival times), see 
\cite[Cor. 2.1]{bernackaite:siaulys:2015}. In such cases, (namely, when \eqref{eq.CKP29V.3.2} holds for any fixed $T>0$), and following the path of our proof, we obtain \eqref{eq.CKP29V.3.3} local uniformly for $T \in \Lambda :=\{t\;:\;\PP(\Theta_1 \leq t\,,\;\Delta_1 \leq t) >0 \}$ (that means uniformly for $T \in \Lambda_{T^*}:= [0,\,T^*]\cap \Lambda$, for any fixed $T^* \in \Lambda$).
\ere

The following result provides a more elegant asymptotic relation than that of \eqref{eq.CKP29V.3.3} when the distributions $F,\,G$ are restricted to the class of regular variation. Its proof is following immediately from the application of dominated convergence theorem and relation \eqref{eq.CKP29V.3.3}, and for sake of brevity is omitted.

\bco \label{cor.CKP29V.3.1}
Let consider the discounted aggregate claims from \eqref{eq.CKP29V.1.2}. We suppose the conditions of Theorem \ref{th.CKP29V.3.1} are valid, with $T>0$, with the restrictions $F \in \mathcal{R}_{-\alpha}$, $G \in \mathcal{R}_{-\beta}$, with $\alpha,\,\beta \geq 0$. Then it holds
\beao
\PP[D_1(T) > x\,,\;D_2(T) >y] \sim \bF(x)\,\bG(y)\,\int_0^T \int_0^T e^{-\alpha\,r_1\,s}\,e^{-\beta\,r_2\,t}\,\E[N(ds)\,M(dt)]\,.
\eeao
\eco

For the case of infinite horizon, we need strictly positive interest force on both lines of business, in order to be satisfied the convergence conditions for the infinite sums of the discounted aggregate claims. Additionally, we should use the following condition, related with the convergence of the counting processes, which depends in some sense on the Matuszewska indexes of $F$ and $G$. 

\begin{assumption} \label{ass.CKP29V.3.2}
Let $r_1 \wedge r_2 >0$ and $F,\,G \in \mathcal{C} \cap \mathcal{P_D}$. We suppose that there exist pairs $(p_1,\,p_2)$ and $(q_1,\,q_2)$, such that 
$0<p_1<J_F^- \leq J_F^+ < p_2 < \infty$,  $0<q_1<J_G^- \leq J_G^+ < q_2 < \infty$, and some $\gamma > 1$ such that it holds
\beam \label{eq.CKP29V.3.4} 
\sum_{i=1}^{\infty} i^{\gamma\,p_2}\,\E^{1/2}\left[ e^{-2\,p_1\,r_1\,\tau_i}\right] < \infty\,,\;\sum_{j=1}^{\infty} j^{\gamma\,q_2}\,\E^{1/2}\left[ e^{-2\,q_1\,r_2\,\sigma_j}\right] < \infty\,.
\eeam
\end{assumption}

The following example shows that a general family of quasi-renewal processes satisfy the convergence conditions \eqref{eq.CKP29V.3.4}  of Assumption \ref{ass.CKP29V.3.2}, and evenmore for any finite $\gamma >1$. We note that in this example we assume that the other conditions of Assumption \ref{ass.CKP29V.3.2} holds, and we prove the validation of \eqref{eq.CKP29V.3.4}.

\bexam \label{exam.CKP29V.a} 
Let consider that $\{\Theta_i\,,\;i \in \bbn \}$ are identically distributed and widely lower orthant dependent (WLOD) (see \cite{wang:wang:gao:2013} for introduction and examples of such kind of dependence), that means there exists a sequence of positive numbers $\{g_L^{\Theta}(n)\,,\;n \in \bbn\}$, called dominating sequence, such that for any 
$n \in \bbn$, $x_1,\,\ldots,\,x_n \in \bbr$ it holds
\beam \label{eq.CKP29V.A}
\PP\left(\bigcap_{i=1}^n \{\Theta_i \leq x_i\} \right) \leq g_L^{\Theta}(n)\,\prod_{i=1}^n \PP(\Theta_i \leq x_i)\,,
\eeam
and further we suppose that there exists some $0 < \vep < -\log \E^{1/2}[e^{-2\,r_1\,p_1\,\Theta_1 }]$, such that it holds
\beam \label{eq.CKP29V.B}
\limsup_{n \to \infty} g_L^{\Theta}(n)\,e^{-\vep\,n} < \infty\,.
\eeam

Similarly, we suppose that $\{\Delta_j\,,\;j \in \bbn \}$ are identically distributed and WLOD with dominating sequence $\{g_L^{\Delta}(n)\,,\;n \in \bbn\}$, 
and there exists some $0 < \vep^*< -\log \E^{1/2}[e^{-2\,r_2\,q_1\,\Delta_1}]$ such that it holds
\beao
\limsup_{n \to \infty} g_L^{\Delta}(n)\,e^{-\vep^*\,n} < \infty\,.
\eeao
Then the convergences in \eqref{eq.CKP29V.3.4} are true.

\pr~
We shall prove only the convergence of the first sum, since the other follows similarly. Due to the fact that the $\{\Theta_i\,,\; i \in \bbn\}$ are 
identically distributed, we obtain that for any fixed $\gamma >1$ it holds
\beao
&&\sum_{i=1}^{\infty} i^{\gamma\,p_2}\,\E^{1/2}[e^{-2\,p_1\,r_1\,\tau_i}] = \sum_{i=1}^{\infty} i^{\gamma\,p_2}\,\E[e^{-2\,p_1\,r_1\,\Theta_1\,i}] \\[2mm]
&& \leq \sum_{i=1}^{\infty} i^{\gamma\,p_2}\,g_L(i)\,\left(\E^{1/2}[e^{-2\,p_1\,r_1\,\Theta_1}] \right)^i < \infty\,,
\eeao
where at the second step we used \cite[Prop. 1.1(1)]{wang:wang:gao:2013}. The convergence at the last sum follows by relation \eqref{eq.CKP29V.B} in combination with the fact that 
$0 < \vep < -\log \E^{1/2}[e^{-2\,p_1\,r_1\,\theta_1}]$.
~\halmos
\eexam

The following theorem, extends \cite[Lem. 5.3]{yang:li:2017}, with respect to (1) - (3) of Subsection \ref{subsec.CKP29V.1.2}.

\bth \label{th.CKP29V.3.2}
Let consider the discounted aggregate claims of relation \eqref{eq.CKP29V.1.2}. We suppose that Assumptions \ref{ass.CKP29V.1.1}, \ref{ass.CKP29V.3.2} are valid and each of the sequences of r.v.s $\{X_i,\,i \in \bbn\}$ and $\{Y_j,\,j \in \bbn\}$  satisfies Assumption \ref{ass.CKP29V.2.2}. Then it holds 
\beam \label{eq.CKP29V.3.5}
\PP[D_1(\infty) > x\,,\;D_2(\infty) >y] \sim \int_0^{\infty} \int_0^{\infty} \bF(x\,e^{r\,s})\,\bG(y\,e^{r\,t})\,\E[N(ds)\,M(dt)]\,.
\eeam   
\ethe

The following corollary provides a more explicit expression than that of \eqref{eq.CKP29V.3.5}, in the case when the distributions $F,\,G$ are regularly varying.

\bco \label{cor.CKP29V.3.2}
Let consider the discounted aggregate claims from \eqref{eq.CKP29V.1.2}. We suppose that the conditions of Theorem \ref{th.CKP29V.3.2} are true, with the restrictions $F \in\mathcal{R}_{-\alpha}$, $G \in\mathcal{R}_{-\beta}$, with $\alpha \wedge \beta >0$. Then it holds
\beam \label{eq.CKP29V.3.6} \notag
&&\PP[D_1(\infty) > x\,,\;D_2(\infty) >y] \sim \bF(x)\,\bG(y) \sum_{i=1}^{\infty} \sum_{j=1}^{\infty} \E[e^{-\alpha\,r_1\,\tau_i}\,e^{-\beta\,r_2\,\sigma_j}]\\[2mm]
&&= \bF(x)\,\bG(y) \int_0^{\infty} \int_0^{\infty} e^{-\alpha\,r_1\,s}\,e^{-\beta\,r_2\,t}\,\E[N(ds)\,M(dt)]\,.
\eeam
\eco

\section{Argumentation} \label{sec.CKP29V.4}

\subsection{Finite horizon} \label{subsec.CKP29V.4.1}

Before proving Theorem \ref{th.CKP29V.3.1}, we need two preliminary lemmas. The first one can be found in \cite[Lem. 2.1]{yang:liu:huang:ma:2012}.

\ble \label{lem.CKP29V.3.1}
Let $\{Z_i\,,\;i \in \bbn\}$ represent non-negative, identical random variables with common distribution $V \in \mathcal{S}$, that satisfy Assumption \ref{ass.CKP29V.2.1}. 
Then for any $\vep > 0$, there exists some constant $K=K(\vep) >0$ such that 
\beam \label{eq.CKP29V.2.6b}
\PP\left(\sum_{i=1}^n Z_i >x \right) \leq K\,(1+\vep)^n\,\bV(x)\,,
\eeam
for any $x \geq 0$ and $n \in \bbn$. 
\ele

The following lemma was given by \cite[Lem. 4.1]{chen:li:cheng:2023} and it uses the following version of dependence, from  Assumption \ref{ass.CKP29V.2.1}, for finite number of random variables, which was introduced by \cite{ko:tang:2008}.

\begin{assumption} \label{ass.CKP29V.4.A^*}
Let $Z_1,\,\ldots,\,Z_n$, $n\in\bbn$ r.v.s. We suppose that there exists constants $C>0$, $x_0 >0$, such that for any $2 \leq l \leq n$ holds
\beao
\sup_{x \geq x_0} \sup_{a \in [x_0,\,x]} \dfrac{\PP\left[\sum_{i=1}^{l-1} Z_i > x-a\;\big|\;Z_l = a \right]}{\PP\left[\sum_{i=1}^{l-1} Z_i > x-a \right]} < C\,.
\eeao
\end{assumption}

\ble \label{lem.CKP29V.3.2}
Let $Z_1,\,\ldots,\,Z_n\,,\;n \in \bbn$, non-negative, random variables with distributions $V_1,\,\ldots,\,V_n \in \mathcal{L}$, that satisfy Assumption \ref{ass.CKP29V.4.A^*}. If there exists some distribution $V \in \mathcal{S}$, such that $\bV_i(x) \asymp \bV(x)$, as $\xto$, for any $1\leq i \leq n$, then for any $0<a\leq b <\infty$, it holds
\beao
\PP\left(\sum_{i=1}^n c_i\,Z_i >x  \right) \sim \sum_{i=1}^n\PP\left( c_i\,Z_i >x  \right)\,,
\eeao
as $\xto$, uniformly for $(c_1,\,\ldots,\,c_n) \in [a,\,b]^n$.
\ele

The line of proof of Theorem \ref{th.CKP29V.3.1}, in spite of the much more general assumptions, looks like similar to that of \cite[Lem. 4.4]{yang:li:2017}, with the necessary changes due to our general set up. 

\noindent{\bf Proof of Theorem \ref{th.CKP29V.3.1}.}~
We choose some large enough $N \in \bbn$ and for $x\wedge y >0$ we obtain
\beam \label{eq.CKP29V.2.7b}
&&\PP[D_1(t)  > x\,,\;D_2(t) >y] = \PP\left(\sum_{i=1}^{N(t)} X_i\,e^{-r_1\,\tau_i} >x\,,\;\sum_{j=1}^{M(t)} Y_j\,e^{-r_2\,\sigma_j} >y\right) \\[2mm] \notag
&&=\left[\left(\sum_{n=1}^{\infty} \sum_{m=1}^{\infty} -\sum_{n=1}^{N} \sum_{m=1}^{N}  \right) + \sum_{n=1}^{N} \sum_{m=1}^{N} \right] \PP\Bigg(\sum_{i=1}^{n} X_i\,e^{-r_1\,\tau_i} >x\,,\;\sum_{j=1}^{m} Y_j\,e^{-r_2\,\sigma_j} >y\,, \\[2mm] \notag
&&\;N(T)=n\,,\; M(T)=m\Bigg)=J_1(x,\,T,\,N) + J_2(x,\,T,\,N)\,.
\eeam 

Let start with $J_1(x,\,T,\,N)$. Initially we obtain that it holds
\beam \label{eq.CKP29V.2.8}
&&J_1(x,\,T,\,N) \leq\left(\sum_{n=1}^{\infty} \sum_{m=N+1}^{\infty} +\sum_{n=N+1}^{\infty} \sum_{m=1}^{\infty}  \right)  \\[2mm] \notag
&&\PP\left(\sum_{i=1}^{n} X_i\,e^{-r_1\,\tau_i} >x\,,\;\sum_{j=1}^{m} Y_j\,e^{-r_2\,\sigma_j} >y\,,\;N(T)=n\,,\; M(T)=m\right) \\[2mm] \notag
&&=J_{11}(x,\,T,\,N)+J_{12}(x,\,T,\,N)\,.
\eeam
For $J_{11}(x,\,T,\,N)$, for any $\vep >0$, we find
\beam \label{eq.CKP29V.2.9} \notag
&&J_{11}(x,\,T,\,N) \leq \sum_{n=1}^{\infty} \sum_{m=N+1}^{\infty} \PP\left(\sum_{i=1}^{n} X_i\,e^{-r_1\,\tau_1} >x\,,\;\sum_{j=1}^{m} Y_j\,e^{-r_2\,\sigma_1} >y\,,\;\tau_n \leq T\,,\; \sigma_m \leq T\right) \\[2mm] 
&&=  \sum_{n=1}^{\infty} \sum_{m=N+1}^{\infty} \int_0^T \int_0^T \PP\left(\sum_{i=1}^{n} X_i\,e^{-r_1\,s} >x\,,\;\sum_{j=1}^{m} Y_j\,e^{-r_2\,t} >y\right)\\[2mm] \notag
&&\times \,\PP\left(\sum_{i=2}^{n} \Theta_i \leq T-s\,,\;\sum_{j=2}^{m} \Delta_j \leq T-t \right)\,\PP(\Theta_1 \in ds\,,\;\Delta_1 \in dt) \\[2mm] \notag
&&\leq \sum_{n=1}^{\infty} \sum_{m=N+1}^{\infty} \int_0^T \int_0^T \PP\left(\sum_{i=1}^{n} X_i >x\,e^{r_1\,s}\right)\,\PP\left(\sum_{j=1}^{m} Y_j >y\,e^{r_2\,t}\right)\\[2mm] \notag
&&\times \,\PP\left(N^*(T) \geq n-1\,,\;M^*(T) \geq m-1 \right)\,\PP(\Theta_1 \in ds\,,\;\Delta_1 \in dt) \\[2mm] \notag
&&\leq K_1\,K_2\,\sum_{n=1}^{\infty} \sum_{m=N+1}^{\infty} (1+\vep)^{n+m}\,\PP\left(N^*(T) \geq n-1\,,\;M^*(T) \geq m-1 \right)\\[2mm] \notag
&&\times \,\int_0^T \int_0^T \bF\left(x\,e^{r_1\,s}\right)\,\bG\left(y\,e^{r_2\,t}\right)\,\PP(\Theta_1 \in ds\,,\;\Delta_1 \in dt) \\[2mm] \notag
&&\leq K_1\,K_2\,\left(\sum_{n=1}^{\infty} \sum_{m=N+1}^{\infty} (1+\vep)^{n+m}\,\PP\left(N^*(T) \geq n-1\,,\;M^*(T) \geq m-1 \right) \right)\\[2mm] \notag
&&\times \,\int_0^T \int_0^T \bF\left(x\,e^{r_1\,s}\right)\,\bG\left(y\,e^{r_2\,t}\right)\,\E\left[N(ds)\,M(dt)\right]\,,
\eeam
where at the third step we used Assumption \ref{ass.CKP29V.1.1} and the definition of delayed counting process from relation \eqref{eq.CKP29V.3.1}. At the fourth step the constants $K_1=K_1(\vep) >0$ and $K_2=K_2(\vep) >0$, are given by the application of Lemma \ref{lem.CKP29V.3.1} twice. Further we have
\beao
&&\sum_{n=1}^{\infty} \sum_{m=N+1}^{\infty} (1+\vep)^{n+m}\,\PP\left(N^*(T) \geq n-1\,,\;M^*(T) \geq m-1 \right) \\[2mm]
&&=\sum_{n=1}^{\infty} \sum_{m=N+1}^{\infty} (1+\vep)^{n+m}\,\sum_{i=n-1}^{\infty} \sum_{j=m-1}^{\infty} \PP\left(N^*(T) =i\,,\;M^*(T) =j \right)  \\[2mm]
&&=\sum_{i=0}^{\infty} \sum_{j=N}^{\infty} \, \PP\left(N^*(T) =i\,,\;M^*(T) =j \right)\,\sum_{n=1}^{i+1}(1+ \vep)^n \sum_{m=N+1}^{j+1} (1+\vep)^{m}  \\[2mm]
&&\leq \sum_{i=0}^{\infty} \sum_{j=N}^{\infty} \dfrac{(1+ \vep)^{i+2} (1+\vep)^{j+2}}{\vep^2}\,\PP\left(N^*(T) =i\,,\;M^*(T) =j \right) \\[2mm]
&&= \dfrac{(1+ \vep)^4}{\vep^2}\, \E \left[ (1+\vep)^{N^*(T)}\,(1+\vep)^{M^*(T)}\,{\bf 1}_{\{M^*(T) \geq N\}}\right] \\[2mm]
&&= \dfrac{(1+ \vep)^4}{\vep^2}\, \left(\E^{1/2} \left[ (1+\vep)^{2\,N^*(T)}\right]\,\E^{1/2} \left[ (1+\vep)^{2\,M^*(T)}\,{\bf 1}_{\{M^*(T) \geq N\}}\right] \right)\,,
\eeao
where at the second step we changed the order of summation, while at the last step we applied the H\"{o}lder inequality. From relation \eqref{eq.CKP29V.3.2}, we can find a properly chosen $\vep > 0$ such that it holds
\beao
\E\left[ (1+\vep)^{2\,N^*(T)}\right] \leq \E \left[ e^{\beta\,N^*(T)}\right] < \infty\,,
\eeao
and similarly for the $\{M^*(T)\,,\; T \geq 0\}$. Hence, in combination with the last relation and \eqref{eq.CKP29V.2.9}, is implied that for any $\delta \in (0,\,1)$, there exists some large enough $N \in \bbn$, such that 
\beam \label{eq.CKP29V.2.10}
J_{11}(x,\,T,\,N) \leq \dfrac {\delta}2 \int_0^T \int_0^T \bF\left(x\,e^{r_1\,s} \right) \,\bG\left(y\,e^{r_2\,t}\right)\,\E\left[N(ds)\,M(dt)\right] \,,
\eeam
and with symmetric argumentation, it holds
\beam \label{eq.CKP29V.2.11}
J_{12}(x,\,T,\,N) \leq \dfrac {\delta}2 \int_0^T \int_0^T \bF\left(x\,e^{r_1\,s} \right) \,\bG\left(y\,e^{r_2\,t}\right)\,\E\left[N(ds)\,M(dt)\right] \,.
\eeam

From \eqref{eq.CKP29V.2.10} and \eqref{eq.CKP29V.2.11}, in combination with \eqref{eq.CKP29V.2.8}, we obtain that
\beam \label{eq.CKP29V.2.12}
J_{1}(x,\,T,\,N) \leq \delta \int_0^T \int_0^T \bF\left(x\,e^{r_1\,s} \right) \,\bG\left(y\,e^{r_2\,t}\right)\,\E\left[N(ds)\,M(dt)\right] \,.
\eeam

Before proceed to deal with $J_{2}(x,\,T,\,N)$, we note that for any fixed $n,\,m \in \bbn$, conditioning with respect to values of $\tau_1,\,\ldots,\,\tau_{n+1},\,\sigma_1,\,\ldots,\,\sigma_{m+1}$, we get
\beam \label{eq.CKP29V.2.13}
&& \PP\left(\sum_{i=1}^{n} X_i\,e^{-r_1\,\tau_i} >x\,,\;\sum_{j=1}^{m} Y_j\,e^{-r_2\,\sigma_j} >y\,,\;N(T) =n\,,\; M(T)= m\right) \\[2mm]  \notag
&&= \int\cdots\int_{0\leq s_1 \leq \cdots \leq s_n \leq T < s_{n+1},\; 0\leq t_1 \leq \cdots \leq t_m \leq T <t_{m+1}} \PP\left(\sum_{i=1}^{n} X_i\,e^{-r_1\,s_i} >x\right) \\[2mm] \notag
&&\times \,\PP\left(\sum_{j=1}^{m} Y_j\,e^{-r_2\,t_j} >y\right)\,\PP\left(\tau_1 \in ds_1,\,\ldots,\,\tau_{n+1}\in ds_{n+1},\,\sigma_1 \in dt_1,\,\ldots,\,\sigma_{m+1}\in ds_{m+1}\right) \\[2mm] \notag
&&\sim \sum_{i=1}^{n} \sum_{j=1}^{m} \int\cdots\int_{0\leq s_1 \leq \cdots \leq s_n \leq T < s_{n+1},\; 0\leq t_1 \leq \cdots \leq t_m \leq T <t_{m+1}}\,\bF\left(x\,e^{r_1\,s_i}\right)\,\bG\left(y\,e^{r_2\,t_i}\right)\\[2mm] \notag
&&\times \,\PP\left(\tau_1 \in ds_1,\,\ldots,\,\tau_{n+1} \in ds_{n+1},\,\sigma_1 \in dt_1,\,\ldots,\,\sigma_{m+1} \in dt_{m+1}\right) \\[2mm] \notag
&& =\,\sum_{i=1}^{n} \sum_{j=1}^{m}\,\PP\left(X_i\,e^{-r_1\,\tau_i} >x\,,\;Y_j\,e^{-r_2\,\sigma_j} >y\,,\;N(T) =n\,,\; M(T)= m\right)\,,
\eeam
where at the second step we applied Lemma \ref{lem.CKP29V.3.2} twice. From relation \eqref{eq.CKP29V.2.13} we obtain
\beam \label{eq.CKP29V.2.14}
&&J_{2}(x,\,T,\,N) \\[2mm]  \notag
&&\sim \sum_{n=1}^{N} \sum_{m=1}^{N} \sum_{i=1}^{n} \sum_{j=1}^{m}
\,\PP\left(X_i\,e^{-r_1\,\tau_i} >x\,,\; Y_j\,e^{-r_2\,\sigma_j} >y\,,\;N(T) =n\,,\; M(T)= m\right) \\[2mm] \notag
&&= \left[ \sum_{n=1}^{\infty} \sum_{m=1}^{\infty} - \left(\sum_{n=1}^{\infty} \sum_{m=1}^{\infty} -\sum_{n=1}^{N} \sum_{m=1}^{N} \right) \right] \sum_{i=1}^{n} \sum_{j=1}^{m} \,\PP\big( X_i\,e^{-r_1\,\tau_i} >x\,,\\[2mm]  \notag
&&\; Y_j\,e^{-r_2\,\sigma_j} >y\,,\;N(T) =n\,,\; M(T)= m\big) =:J_{21}(x,\,T,\,N) - J_{22}(x,\,T,\,N)\,.
\eeam
Changing the order of summation in $J_{21}(x,\,T,\,N)$ we find
\beam \label{eq.CKP29V.2.15}
&&J_{21}(x,\,T,\,N) \\[2mm]  \notag
&&= \sum_{i=1}^{\infty} \sum_{j=1}^{\infty}\,\PP\left( X_i\,e^{-r_1\,\tau_i} >x\,,\; Y_j\,e^{-r_2\,\sigma_j} >y\,,\;\tau_i \leq T\,,\; \sigma_j \leq T\right) \\[2mm] \notag
&&= \sum_{i=1}^{\infty} \sum_{j=1}^{\infty} \int_0^T \,\int_0^T \,\PP\big( X_i\,e^{-r_1\,s} >x\,,\; Y_j\,e^{-r_2\,t} >y\big)\,\PP(\tau_i \in ds\,,\;\sigma_j \in dt)\\[2mm] \notag
&&= \int_0^T \,\int_0^T \,\bF\big(x\,e^{r_1\,s}\big)\,\,\bG\big(y\,e^{r_2\,t}\big)\,\E\left[N(ds)\,M(dt)\right] \,.
\eeam

Now we deal with $J_{22}(x,\,T,\,N)$. It holds
\beam \label{eq.CKP29V.2.16}
&&J_{22}(x,\,T,\,N)\leq \left(\sum_{n=1}^{\infty} \sum_{m=N+1}^{\infty} + \sum_{n=N+1}^{\infty} \sum_{m=1}^{\infty} \right)\, \\[2mm]  \notag
&&\sum_{i=1}^{n} \sum_{j=1}^{m}\PP\left( X_i\,e^{-r_1\,\tau_i} >x\,,\; Y_j\,e^{-r_2\,\sigma_j} >y\,,\;\tau_n \leq T\,,\; \sigma_m \leq T\right) \\[2mm] \notag
&&= J_{221}(x,\,T,\,N) + J_{222}(x,\,T,\,N)\,.
\eeam
For $J_{221}(x,\,T,\,N)$ we obtain
\beam \label{eq.CKP29V.2.17}
&&J_{221}(x,\,T,\,N) =\sum_{n=1}^{\infty} \sum_{m=N+1}^{\infty} \sum_{i=1}^{n} \sum_{j=1}^{m}\\[2mm]  \notag
&&\,\int_0^T \int_0^T \,\bF\big(x\,e^{r_1\,s}\big)\,\,\bG\big(y\,e^{r_2\,t}\big)\,\PP\left(N^*(T-s) \geq n-1\,,\;M^*(T -t) \geq m-1\right)\, \\[2mm]  \notag
&&\times \PP(\Theta_1 \in ds\,,\;\Delta_1 \in dt)\leq \left(\sum_{n=1}^{\infty} \sum_{m=N+1}^{\infty} n\,m \,\PP\left(N^*(T) \geq n-1\,,\;M^*(T) \geq m-1\right) \,\right)  \\[2mm]  \notag
&&\times \int_0^T \int_0^T \,\bF\big(x\,e^{r_1\,s}\big)\,\,\bG\big(y\,e^{r_2\,t}\big)\,\E\left[N(ds)\,M(dt)\right] \,.
\eeam

Further we have 
\beam \label{eq.CKP29V.2.18}
&&\sum_{n=1}^{\infty} \sum_{m=N+1}^{\infty} n\,m \,\PP\left(N^*(T) \geq n-1\,,\;M^*(T) \geq m-1\right)  \\[2mm]  \notag
&&= \sum_{n=1}^{\infty} \sum_{m=N+1}^{\infty} n\,m \,\sum_{i=n-1}^{\infty} \sum_{j=m-1}^{\infty} \PP\left(N^*(T) =i\,,\;M^*(T) =j\right)  \\[2mm]  \notag
&&=\,\sum_{i=0}^{\infty} \sum_{j=N}^{\infty} \PP\left(N^*(T) =i\,,\;M^*(T) =j\right) \sum_{n=1}^{i+1}  n\,\sum_{m=N+1}^{j+1} m   \\[2mm]  \notag
&& \leq \sum_{i=0}^{\infty} \sum_{j=N}^{\infty} (1+ i)^2\,(1+j)^2\, \PP\left(N^*(T) =i\,,\;M^*(T) =j\right)  \\[2mm]  \notag
&& =\E\left[(N^*(T) +1)^2\,(M^*(T)+1)^2\,{\bf 1}_{\{M^*(T) \geq N\}}\right]  \\[2mm]  \notag
&&\leq \E^{1/2}\left[(N^*(T) +1)^4\right] \,\E^{1/2}\left[(M^*(T)+1)^4\,{\bf 1}_{\{M^*(T) \geq N\}}\right] \,,
\eeam
where at the second step we changed the order of summation, while at the last step we applied the H\"{o}lder inequality. From the Assumption \eqref{eq.CKP29V.3.2}, for any $\delta \in (0,\,1)$ we can find some large enough $N$, such that from relations \eqref{eq.CKP29V.2.17} and \eqref{eq.CKP29V.2.18} is implied
\beam \label{eq.CKP29V.2.19}
J_{221}(x,\,T,\,N) \leq \dfrac {\delta}2  \int_0^T \int_0^T \,\bF\big(x\,e^{r_1\,s}\big)\,\,\bG\big(y\,e^{r_2\,t}\big)\,\E\left[N(ds)\,M(dt)\right]\,.
\eeam

Similarly, due to symmetry, it holds
\beam \label{eq.CKP29V.2.20}
J_{222}(x,\,T,\,N) \leq \dfrac {\delta}2  \int_0^T \int_0^T \,\bF\big(x\,e^{r_1\,s}\big)\,\,\bG\big(y\,e^{r_2\,t}\big)\,\E\left[N(ds)\,M(dt)\right]\,.
\eeam
From relations \eqref{eq.CKP29V.2.19} and \eqref{eq.CKP29V.2.20}, in combination with relation \eqref{eq.CKP29V.2.16} we obtain
\beam \label{eq.CKP29V.2.21}
J_{22}(x,\,T,\,N) \leq \delta \int_0^T \int_0^T \,\bF\big(x\,e^{r_1\,s}\big)\,\,\bG\big(y\,e^{r_2\,t}\big)\,\E\left[N(ds)\,M(dt)\right]\,.
\eeam
Putting \eqref{eq.CKP29V.2.15} and \eqref{eq.CKP29V.2.21} into \eqref{eq.CKP29V.2.14} we find
\beam \label{eq.CKP29V.2.22}
&&(1-\delta) \int_0^T \int_0^T \,\bF\big(x\,e^{r_1\,s}\big)\,\,\bG\big(y\,e^{r_2\,t}\big)\,\E\left[N(ds)\,M(dt)\right] \lesssim\\[2mm] \notag
&&J_{2}(x,\,T,\,N) \lesssim \int_0^T \int_0^T \,\bF\big(x\,e^{r_1\,s}\big)\,\,\bG\big(y\,e^{r_2\,t}\big)\,\E\left[N(ds)\,M(dt)\right]\,,
\eeam
Putting \eqref{eq.CKP29V.2.12} and \eqref{eq.CKP29V.2.22} into \eqref{eq.CKP29V.2.7b} we obtain
\beao
&&(1-\delta)\,\int_0^T \int_0^T \,\bF\big(x\,e^{r_1\,s}\big)\,\,\bG\big(y\,e^{r_2\,t}\big)\,\E\left[N(ds)\,M(dt)\right] \lesssim\\[2mm] \notag
&&\PP[D_1(T)  > x\,,\;D_2(T) >y] \lesssim (1+\delta)\,\int_0^T \int_0^T \,\bF\big(x\,e^{r_1\,s}\big)\,\,\bG\big(y\,e^{r_2\,t}\big)\,\E\left[N(ds)\,M(dt)\right]\,,
\eeao
from which by letting $\delta \downarrow 0$ gives relation \eqref{eq.CKP29V.3.3}. 
~\halmos

\subsection{Infinite horizon} \label{subsec.CKP29V.4.2}

Before the proofs of Theorem \ref{th.CKP29V.3.2} and Corollary \ref{cor.CKP29V.3.2}, we need some preliminary lemmas. For sake of brevity, hereafter we shall use the $(p_1,\,p_2)\,,\;(q_1,\,q_2)\,,\;\gamma$ as they were presented in Assumption \ref{ass.CKP29V.3.2}. For all the applications of Potter inequalities 
\eqref{eq.CKP29V.2.3}, \eqref{eq.CKP29V.2.4}, according to $p_1,\,q_1$ and $p_2,\,q_2$ respectively, we shall use common constants $C_1,\,D_1$, for inequality \eqref{eq.CKP29V.2.3} and common constants $C_2,\,D_2$, for inequality \eqref{eq.CKP29V.2.4}, for any $V \in \{F,\,G\}$, without loss of generality (for example in \eqref{eq.CKP29V.2.3} and \eqref{eq.CKP29V.2.4} for $F,\,G$, it is enough to choose $D_1=D_1^{(F)} \vee D_1^{(G)}$, $D_2=D_2^{(F)} \vee D_2^{(G)}$, $C_1=C_1^{(F)} \wedge C_1^{(G)}$ and $C_2=C_2^{(F)} \vee C_2^{(G)}$, where the $D_i^{(F)}$, $C_i^{(F)}$, with $i=1,\,2$ correspond to distribution $F$ and the $D_i^{(G)}$, $C_i^{(G)}$, with $i=1,\,2$ correspond to distribution $G$). 
  
\ble \label{lem.CKP29V.4.3} 
Under the conditions of Theorem \ref{th.CKP29V.3.2}, with the more relaxed conditions $F,\,G \in \mathcal{D}\cap \mathcal{A}$, for all $i,\,j \in \bbn$ it holds
\beam \label{eq.CKP29V.4.18} 
\PP(X_i\,e^{-r_1\,\tau_i} > x\,,\;Y_j\,e^{-r_2\,\sigma_j} >y) \asymp \bF(x)\,\bG(y)\,.
\eeam
\ele

\pr~
At first we shall show that
\beam \label{eq.CKP29V.4.19}
\limsup \dfrac{\PP(X_i\,e^{-r_1\,\tau_i} > x\,,\;Y_j\,e^{-r_2\,\sigma_j} >y)}{\bF(x)\,\bG(y)} < \infty\,.
\eeam
Since $F,\,G \in \mathcal{D}\cap \mathcal{A} \subsetneq \mathcal{D} \cap \mathcal{P_D}$, from relation \eqref{eq.CKP29V.2.3}, for any $x \geq D_1$ and $i,\,j \in \bbn$ we obtain 
\beao
&&\PP(X_i\,e^{-r_1\,\tau_i} > x\,,\;Y_j\,e^{-r_2\,\sigma_j} >y)=\int_0^{\infty} \int_0^{\infty} 
\bF(x\,e^{r_1\,s})\,\bG(y\,e^{r_2\,t})\,\PP(\tau_i \in ds\,,\;\sigma_j \in dt) \\[2mm]
&&\leq \dfrac 1{C_1^2}\,\bF(x)\,\bG(y)\,\int_0^{\infty} \int_0^{\infty} 
e^{-r_1\,p_1\,s}\,e^{-r_2\,q_1\,t}\,\PP(\tau_i \in ds\,,\;\sigma_j \in dt) \\[2mm]
&&= \dfrac 1{C_1^2}\,\bF(x)\,\bG(y)\,\E[e^{-r_1\,p_1\,\tau_i}\,e^{-r_2\,q_1\,\sigma_j}] \leq \dfrac 1{C_1^2}\,\bF(x)\,\bG(y) \,,
\eeao
where in last step, we recall that $r_1 \wedge r_2 >0$, that implies \eqref{eq.CKP29V.4.19}.

From the other hand side, by \eqref{eq.CKP29V.2.4} follows that for any $x \geq D_2$ it holds
\beao
&&\PP(X_i\,e^{-r_1\,\tau_i} > x\,,\;Y_j\,e^{-r_2\,\sigma_j} >y) \\[2mm]
&&=\int_0^{\infty} \int_0^{\infty} 
\bF(x\,e^{r_1\,s})\,\bG(y\,e^{r_2\,t})\,\PP(\tau_i \in ds\,,\;\sigma_j \in dt) \\[2mm]
&&\geq \dfrac 1{C_2^2}\,\bF(x)\,\bG(y)\,\int_0^{\infty} \int_0^{\infty} 
e^{-r_1\,p_2\,s}\,e^{-r_2\,q_2\,t}\,\PP(\tau_i \in ds\,,\;\sigma_j \in dt) \\[2mm]
&&= \dfrac 1{C_2^2}\,\bF(x)\,\bG(y)\,\E[e^{-r_1\,p_2\,\tau_i}\,e^{-r_2\,q_2\,\sigma_j}]  \,,
\eeao
because $\Theta_i,\,\Delta_j$ follow non-defective distributions, the last expectation is strictly positive. Thus, it holds
\beam \label{eq.CKP29V.4.20}
\liminf \dfrac{\PP(X_i\,e^{-r_1\,\tau_i} > x\,,\;Y_j\,e^{-r_2\,\sigma_j} >y)}{\bF(x)\,\bG(y)} > 0\,,
\eeam
for any $i,\,j \in \bbn$. From \eqref{eq.CKP29V.4.19} and \eqref{eq.CKP29V.4.20} we obtain \eqref{eq.CKP29V.4.18}.
~\halmos

The next lemma plays crucial role in the proof of Lemma \ref{lem.CKP29V.4.5*}, that is important for the proof of Theorem \ref{th.CKP29V.3.2}. In some sense, through this lemma we recognize the necessity of class $\mathcal{P_D}$ in the conditions of the main result on the infinite horizon. 
  
\ble \label{lem.CKP29V.4.4} 
Under the conditions of Theorem \ref{th.CKP29V.3.2}, with the more relaxed conditions $F,\,G \in \mathcal{D}\cap \mathcal{A}$, it holds
\beam \label{eq.CKP29V.4.21} 
&&\lim_{N \to \infty} \limsup \dfrac{\PP\left( \sum_{i=N}^{\infty} X_i\,e^{-r_1\,\tau_i} > x\,,\;\sum_{j=1}^{\infty} Y_j\,e^{-r_2\,\sigma_j} >y\right) }{\bF(x)\,\bG(y)} \\[2mm] \notag
&&= \lim_{N \to \infty} \limsup \dfrac{ \sum_{i=N}^{\infty}\sum_{j=1}^{\infty} \PP\left( X_i\,e^{-r_1\,\tau_i} > x\,,\;Y_j\,e^{-r_2\,\sigma_j} >y\right) }{\bF(x)\,\bG(y)} =0\,,
\eeam
and 
\beam \label{eq.CKP29V.4.22} 
&&\lim_{N \to \infty} \limsup \dfrac{\PP\left( \sum_{i=1}^{\infty} X_i\,e^{-r_1\,\tau_i} > x\,,\;\sum_{j=N}^{\infty} Y_j\,e^{-r_2\,\sigma_j} >y\right) }{\bF(x)\,\bG(y)} \\[2mm] \notag
&&= \lim_{N \to \infty} \limsup \dfrac{ \sum_{i=1}^{\infty}\sum_{j=N}^{\infty} \PP\left( X_i\,e^{-r_1\,\tau_i} > x\,,\;Y_j\,e^{-r_2\,\sigma_j} >y\right) }{\bF(x)\,\bG(y)} =0\,.
\eeam
\ele

\pr~
We shall prove only relation \eqref{eq.CKP29V.4.21}, since  \eqref{eq.CKP29V.4.22} comes similarly via symmetric arguments.

We choose some large enough $K >0$, such that $\sum_{i=1}^{\infty} 1/i^{\gamma} < K$. Then it holds
\beam \label{eq.CKP29V.4.23} 
&&\PP\left( \sum_{i=N}^{\infty} X_i\,e^{-r_1\,\tau_i} > x\,,\;\sum_{j=1}^{\infty} Y_j\,e^{-r_2\,\sigma_j} >y\right)  \\[2mm] \notag
&&\leq \PP\left( \sum_{i=N}^{\infty} 
X_i\,e^{-r_1\,\tau_i} > x\,\sum_{i=N}^{\infty} \dfrac 1{i^{\gamma}\,K}\,,\; \sum_{j=1}^{\infty} Y_j\,e^{-r_2\,\sigma_j} >y\,\sum_{j=1}^{\infty} \dfrac 1{j^{\gamma}\,K}\right) \\[2mm] \notag
&&\leq \PP\left( \bigcup_{i=N}^{\infty} \left\{
X_i\,e^{-r_1\,\tau_i} > x\,\dfrac 1{i^{\gamma}\,K} \right\}\,,\;\bigcup_{j=1}^{\infty} \left\{ Y_j\,e^{-r_2\,\sigma_j} >y\,\dfrac 1{j^{\gamma}\,K} \right\} \right) \\[2mm] \notag
&&\leq \sum_{i=N}^{\infty}\sum_{j=1}^{\infty} \PP\left( X_i\,e^{-r_1\,\tau_i} > x\,\dfrac 1{i^{\gamma}\,K}\,,\;Y_j\,e^{-r_2\,\sigma_j} >y\,\dfrac 1{j^{\gamma}\,K} \right)=: J(x,\,y;\,N) \,.
\eeam

We denote $D:=D_1 \vee D_2$, and we separate $J(x,\,y;\,N)$ as follows:
\beam \label{eq.CKP29V.4.24} \notag 
&&J(x,\,y;\,N)=\sum_{i=N}^{\infty} \sum_{j=1}^{\infty} \int_0^1 \int_0^1 \bF\left(\dfrac x{i^{\gamma}\,K\,v}\right)\,\bG\left(\dfrac y{j^{\gamma}\,K\,u} \right)\,\PP\left(e^{-r_1\,\tau_i} \in dv\,,\;e^{-r_2\,\sigma_j} \in du \right) \\[2mm] \notag
&&=\sum_{i=N}^{\infty} \sum_{j=1}^{\infty} \Bigg( \int_0^{x/i^{\gamma}\,D\,K} \int_0^{y/j^{\gamma}\,D\,K} +\int_0^{x/i^{\gamma}\,D\,K} \int_{y/j^{\gamma}\,D\,K}^1 + \int_{x/i^{\gamma}\,D\,K}^1 \int_0^{y/j^{\gamma}\,D\,K} \\[2mm]
&&+\int_{x/i^{\gamma}\,D\,K}^1 \int_{y/j^{\gamma}\,D\,K}^1 \Bigg) \bF\left(\dfrac x{i^{\gamma}\,K\,v}\right)\,\bG\left(\dfrac y{j^{\gamma}\,K\,u} \right)\,\PP\left(e^{-r_1\,\tau_i} \in dv\,,\;e^{-r_2\,\sigma_j} \in du \right) \\[2mm] \notag
&&=: \sum_{k=1}^4 J_k(x,\,y;\,N)\,.
\eeam
For $J_1(x,\,y;\,N)$, for $x \wedge y \geq D$ and putting $C:=C_2 \vee \dfrac 1{C_1}$ we obtain
\beao
&&J_1(x,\,y;\,N)\leq \sum_{i=N}^{\infty} \sum_{j=1}^{\infty} \int_0^{x/i^{\gamma}\,D\,K} \int_0^{y/j^{\gamma}\,D\,K} C^2\,\bF(x)\,\bG(y)\, \left[ \left(i^{\gamma}\,K\,v\right)^{p_1} + \left(i^{\gamma}\,K\,v\right)^{p_2}\right]\\[2mm]
&&\times \left[\left(j^{\gamma}\,K\,u\right)^{q_1}+\left( j^{\gamma}\,K\,u\right)^{q_2}\right]\,\PP\left(e^{-r_1\,\tau_i} \in dv\,,\;e^{-r_2\,\sigma_j} \in du \right) \\[2mm]
&&\leq \sum_{i=N}^{\infty} \sum_{j=1}^{\infty} \int_0^{1} \int_0^{1} C^2\, \left( 2\,i^{\gamma\,p_2}\,K^{p_2}\,v^{p_1}\right)\,\left( 2\,j^{\gamma\,q_2}\,K^{q_2}\,u^{q_1}\right)\\[2mm]
&&\times \bF(x)\,\bG(y)\,\PP\left(e^{-r_1\,\tau_i} \in dv\,,\;e^{-r_2\,\sigma_j} \in du \right) \\[2mm]
&&=4\,C^2\,K^{p_2+q_2}\,\bF(x)\,\bG(y)\,\sum_{i=N}^{\infty} \sum_{j=1}^{\infty}\,i^{\gamma\,p_2}\,j^{\gamma\,q_2} \E\left[e^{-r_1\,p_1\,\tau_i} \,e^{-r_2\,q_1\,\sigma_j} \right] \\[2mm]
&&\leq 4\,C^2\,K^{p_2+q_2}\,\bF(x)\,\bG(y)\,\sum_{i=N}^{\infty}\,i^{\gamma\,p_2}\, \E^{1/2}\left[e^{-2\,r_1\,p_1\,\tau_i} \right] \sum_{j=1}^{\infty}\,j^{\gamma\,q_2} \E^{1/2}\left[e^{-2\,r_2\,q_1\,\sigma_j} \right]\,,
\eeao
where at the first step we used the inequalities \eqref{eq.CKP29V.2.3}, \eqref{eq.CKP29V.2.4}, while at the last step we applied H\"{o}lder inequality. From Assumption \ref{ass.CKP29V.3.2} we obtain immediately
\beam \label{eq.CKP29V.4.25} 
\lim_{N \to \infty} \limsup \dfrac{J_1(x,\,y;\,N)}{ \bF(x)\,\bG(y)}=0\,.
\eeam
For the second term $J_2(x,\,y;\,N)$, from \eqref{eq.CKP29V.2.3} and \eqref{eq.CKP29V.2.4}, for any $x\geq D$ it holds
\beam \label{eq.CKP29V.4.26} 
&&J_2(x,\,y;\,N)=\sum_{i=N}^{\infty} \sum_{j=1}^{\infty} \int_0^{x/i^{\gamma}\,D\,K} \int_{y/j^{\gamma}\,D\,K}^1 \,\bF\left(\dfrac x{i^{\gamma}\,K\,v}\right)\,\bG\left(\dfrac y{j^{\gamma}\,K\,u} \right)\\[2mm] \notag
&&\times \PP\left(e^{-r_1\,\tau_i} \in dv\,,\;e^{-r_2\,\sigma_j} \in du \right) \leq \sum_{i=N}^{\infty} \sum_{j=1}^{\infty} \int_0^{x/i^{\gamma}\,D\,K} \int_{y/j^{\gamma}\,D\,K}^1 C\,\bF\left(x\right)\,\\[2mm] \notag
&&\left[\left(i^{\gamma}\,K\,v\right)^{p_1}+\left(i^{\gamma}\,K\,v\right)^{p_2}\right]\,C^*\left( \dfrac{u\,j^{\gamma}\,D\,K}{y}\right)^{q_2}\,\PP\left(e^{-r_1\,\tau_i} \in dv\,,\;e^{-r_2\,\sigma_j} \in du \right)\\[2mm] \notag
&& \leq 2\,C\,\,C^*\,K^{p_2+ q_2}\,\left(\dfrac D{y}\right)^{q_2}\,\bF\left(x\right) \sum_{i=N}^{\infty} \sum_{j=1}^{\infty} i^{\gamma\,p_2}\,j^{\gamma\,q_2}\,\E\left[e^{-r_1\,p_1\,\tau_i} \,e^{-r_2\,q_1\,\sigma_j} \right] \\[2mm] \notag
&& \leq 2\,C\,\,C^*\,K^{p_2+ q_2}\,\left(\dfrac D{y}\right)^{q_2}\,\bF\left(x\right) \sum_{i=N}^{\infty}  i^{\gamma\,p_2}\,\E^{1/2}\left[e^{-2\,r_1\,p_1\,\tau_i} \right]\,\sum_{j=1}^{\infty} j^{\gamma\,q_2}\,\E^{1/2}\left[e^{-2\,r_2\,q_1\,\sigma_j} \right] \,,
\eeam
where at the fourth step we applied H\"{o}lder inequality, at the second step, for the distribution $G$, the constant $C^*$ follows by the following argument: Since $G \in \mathcal{D} \cap \mathcal{A}$, there exists some $J_{G}^- \leq \lambda^* \leq J_G^+$, such that we obtain $\E[Y^{\lambda^*}] < \infty$. Hence, there exists $C^* >0$, such that $\E[Y^{\lambda^*}] \leq C^*$. Further, by Markov inequality we find
\beao
\bG\left( \dfrac y{j^{\gamma}\,K\,u} \right) \leq \E\left[Y^{\lambda^*}\right]\,\left( \dfrac{j^{\gamma}\,K\,u}y\right)^{\lambda^*} \leq C^*\,\left( \dfrac{j^{\gamma}\,K\,u\,D}y\right)^{\lambda^*} \leq C^*\,\left( \dfrac{j^{\gamma}\,K\,u\,D}y\right)^{q_2}\,, 
\eeao
where at the last step we took into consideration that $j^{\gamma}\,K\,u\,D \geq y$. Indeed, it is true, since 
\beao
j^{\gamma}\,K\,u\,D \geq \dfrac{j^{\gamma}\,D\,K\,y}{j^{\gamma}\,D\,K} =y\,.
\eeao
From \eqref{eq.CKP29V.4.26}, through Assumption \ref{ass.CKP29V.3.2} and relation \eqref{eq.CKP29V.2.5} we obtain that
\beam \label{eq.CKP29V.4.27} 
\lim_{N \to \infty} \limsup \dfrac{J_2(x,\,y;\,N)}{\bF\left(x\right) \,\bG\left(y\right) }=0 \,.
\eeam

Similarly, for $J_l(x,\,y;\,N)$, with $l=3,\,4$, we can show that
\beam \label{eq.CKP29V.4.28} 
\lim_{N \to \infty} \limsup \dfrac{J_l(x,\,y;\,N)}{\bF\left(x\right) \,\bG\left(y\right) }=0 \,.
\eeam
From relations  \eqref{eq.CKP29V.4.25}, \eqref{eq.CKP29V.4.27}, \eqref{eq.CKP29V.4.28}, in combination with \eqref{eq.CKP29V.4.24} we obtain
\beam \label{eq.CKP29V.4.29} 
\lim_{N \to \infty} \limsup \dfrac{J(x,\,y;\,N)}{\bF\left(x\right) \,\bG\left(y\right) }=0 \,.
\eeam
From  \eqref{eq.CKP29V.4.29} and \eqref{eq.CKP29V.4.23}, we see that the first fraction of  \eqref{eq.CKP29V.4.21} tends to zero. Further, since
\beao
\sum_{i=N}^{\infty} \sum_{j=1}^{\infty} \PP\left( X_i\,e^{-r_1\,\tau_i} > x\,,\;Y_j\,e^{-r_2\,\sigma_j} >y \right)\leq J(x,\,y;\,N)\,,
\eeao
since $1/(j^{\gamma}\,K) <1$, for any $j \in \bbn \cup \{\infty\}$, from \eqref{eq.CKP29V.4.29} is implied that the second fraction of \eqref{eq.CKP29V.4.21} also tends to zero.
~\halmos

The next lemma focuses on multivariate non-linear single big jump principle, for finite randomly weighted sums.

\ble \label{lem.CKP29V.4.5}
Under the conditions of Theorem \ref{th.CKP29V.3.2}, under the more relaxed conditions $F,\,G \in \mathcal{D} \cap \mathcal{L}$, and $r_1,\,r_2 \geq 0$, it holds
\beam \label{eq.CKP29V.4.30} \notag
\PP\left( \sum_{i=1}^{n} X_i\,e^{-r_1\,\tau_i} > x\,,\;\sum_{j=1}^{n}Y_j\,e^{-r_2\,\sigma_j} >y \right) \sim\sum_{i=1}^{n} \sum_{j=1}^{n} \PP\left( X_i\,e^{-r_1\,\tau_i} > x\,,\;Y_j\,e^{-r_2\,\sigma_j} >y \right)\,,\\
\eeam
for any fixed $n \in \bbn$. 
\ele

\pr~
Since the $(X_i,\,Y_j),\,i,\,j \in \bbn$ are independent with $F,\,G \in \mathcal{D} \cap \mathcal{L}$, from \cite[Lem. 3(ii)]{li:2018} we find that there exists a joint insensitivity function $l(\cdot)\;:\;[0,\,\infty) \to (0,\,\infty)$, such that $l(x) < x/2 $, for all $x>0$, $\lim_{\xto} l(x) =\infty$, $l \in \mathcal{R}_0$ and such that for any pair of real numbers $(\kappa,\,\lambda)$ it holds
\beam \label{eq.CKP29V.4.31} \notag
\PP\left( X_i\,e^{-r_1\,\tau_i} > x+\kappa \,l(x)\,,\;Y_j\,e^{-r_2\,\sigma_j} >y + \lambda\,l(y)\right) \sim \PP\left( X_i\,e^{-r_1\,\tau_i} > x\,,\;Y_j\,e^{-r_2\,\sigma_j} >y \right)\,,\\
\eeam
for all $i,\,j \in \bbn$. We note that the necessary moment conditions in \cite[Lem. 3(ii)]{li:2018}, are immediately satisfied for the random weights $e^{-r_1\,\tau_i},\,e^{-r_2\,\sigma_j},\,i,\,j \in \bbn$. We define the events
\beao
A_1=\bigcup_{i=1}^n \left\{X_i\,e^{-r_1\,\tau_i} > x-l(x) \right\}\,,\quad A_2=\bigcup_{j=1}^n \left\{Y_j\,e^{-r_2\,\sigma_j} > y-l(y) \right\}\,,
\eeao
and separate the probability from the left hand side of \eqref{eq.CKP29V.4.30} as follows:
\beam \label{eq.CKP29V.4.32}  \notag
&&\PP\left( \sum_{i=1}^{n} X_i\,e^{-r_1\,\tau_i} > x\,,\;\sum_{j=1}^{n}Y_j\,e^{-r_2\,\sigma_j} >y \right) =\PP\Bigg( \sum_{i=1}^{n} X_i\,e^{-r_1\,\tau_i} > x\,,\;\sum_{j=1}^{n}Y_j\,e^{-r_2\,\sigma_j} >y\,,\\[2mm]
&&\;A_1 \cap A_2 \Bigg)+\PP\left( \sum_{i=1}^{n} X_i\,e^{-r_1\,\tau_i} > x\,,\;\sum_{j=1}^{n}Y_j\,e^{-r_2\,\sigma_j} >y\,,\;A_1^c \cup A_2^c \right)\\[2mm] \notag
&&\leq\PP\Bigg( \sum_{i=1}^{n} X_i\,e^{-r_1\,\tau_i} > x\,,\;\sum_{j=1}^{n}Y_j\,e^{-r_2\,\sigma_j} >y\,,\;A_1 \cap A_2 \Bigg)+ \PP\Bigg( \sum_{i=1}^{n} X_i\,e^{-r_1\,\tau_i} > x\,,\\[2mm] \notag
&& \;\sum_{j=1}^{n}Y_j\,e^{-r_2\,\sigma_j} >y\,,\;A_2^c \Bigg) + \PP\Bigg( \sum_{i=1}^{n} X_i\,e^{-r_1\,\tau_i} > x\,,\;\sum_{j=1}^{n}Y_j\,e^{-r_2\,\sigma_j} >y\,,\;A_1^c \Bigg)\\[2mm] \notag
&&=:\sum_{i=1}^3 L_i(x,\,y,\,n) \,.
\eeam
For the $L_1(x,\,y,\,n)$, through relation \eqref{eq.CKP29V.4.31} we obtain
\beam \label{eq.CKP29V.4.33}  \notag
&&L_1(x,\,y,\,n) \leq \PP\Bigg( A_1 \cap A_2 \Bigg) \leq \sum_{i=1}^{n} \sum_{j=1}^{n} \PP\left(  X_i\,e^{-r_1\,\tau_i} > x-l(x)\,,\;Y_j\,e^{-r_2\,\sigma_j} >y- l(y) \right)\\[2mm] 
&&\sim \sum_{i=1}^{n} \sum_{j=1}^{n} \PP\left( X_i\,e^{-r_1\,\tau_i} > x\,,\;Y_j\,e^{-r_2\,\sigma_j} >y \right) \,.
\eeam
For the $L_2(x,\,y,\,n)$, we have
\beam \label{eq.CKP29V.4.34}  \notag
&&L_2(x,\,y,\,n) = \PP\Bigg( \bigcup_{i=1}^{n} \left\{ X_i\,e^{-r_1\,\tau_i} > \dfrac xn \right\}\,,\; \bigcup_{j=1}^{n} \left\{ Y_j\,e^{-r_2\,\sigma_j} >\dfrac yn \right\}\,,\;\sum_{l=1}^n X_l\,e^{-r_1\,\tau_l} > x\,,\\[2mm] 
&&\;\sum_{k=1}^{n}Y_k\,e^{-r_2\,\sigma_k} >y\,,\; \bigcap_{m=1}^n \left\{Y_m\,e^{-r_2\,\sigma_m}\leq y - l(y) \right\} \Bigg)\\[2mm] \notag
&&\leq \sum_{i=1}^{n} \sum_{j=1}^{n} \PP\Bigg( X_i\,e^{-r_1\,\tau_i} > \dfrac xn\,,\;Y_j\,e^{-r_2\,\sigma_j} > \dfrac yn\,,\;\sum_{k=1}^{n}Y_k\,e^{-r_2\,\sigma_k} >y\,,\\[2mm] \notag
&&\; \bigcap_{m=1}^n \left\{Y_m\,e^{-r_2\,\sigma_m}\leq y - l(y) \right\} \Bigg) \leq \sum_{i=1}^{n} \sum_{j=1}^{n} \sum_{i\neq k=1}^{n} \PP\Bigg( X_i\,e^{-r_1\,\tau_i} > \dfrac xn\,,\;Y_j\,e^{-r_2\,\sigma_j} > \dfrac yn\,,\\[2mm] \notag
&&\;Y_k\,e^{-r_2\,\sigma_k} >\dfrac{l(y)}n \Bigg)
 =o\Bigg[ \sum_{i=1}^{n} \sum_{j=1}^{n} \PP\Bigg( X_i\,e^{-r_1\,\tau_i} >x\,,\;Y_j\,e^{-r_2\,\sigma_j} >y\Bigg) \Bigg]\,,
\eeam
where at the last step we used Assumption \ref{ass.CKP29V.2.2} for the $Y_j,\,j \in \bbn$, together with the inclusions $F,\,G \in \mathcal{D}$. Indeed, for each of the terms of the summation, at the pre-last step of \eqref{eq.CKP29V.4.34}, in combination with Assumption \ref{ass.CKP29V.2.2}, for any $\delta> 0$, there exists a $y_0=y_0(\delta)>0$, and a constant $C^*>0$, due to $F,\,G \in \mathcal{D}$, such that we obtain
\beam \label{eq.CKP29V.4.35}  
&&\PP\Bigg( X_i\,e^{-r_1\,\tau_i} > \dfrac xn \,,\; Y_j\,e^{-r_2\,\sigma_j} >\dfrac yn\,,\;Y_k\,e^{-r_2\,\sigma_k} >\dfrac{l(y)}n  \Bigg)\\[2mm] \notag
&&=\int_0^{\infty} \int_0^{\infty} \int_0^{\infty} \PP\left( X_i\,e^{-r_1\,s} > \dfrac xn \right)\,\PP \left(Y_j\,e^{-r_2\,t_j} > \dfrac yn\,,\;Y_k\,e^{-r_2\,t_k} >\dfrac{l(y)}n \right)\,\\[2mm] \notag
&&\times \PP(\tau_i \in ds\,,\;\sigma_j \in dt_j\,,\; \sigma_k \in dt_k) \leq \delta\, \int_0^{\infty} \int_0^{\infty} \int_0^{\infty} \bF\left(\dfrac{x\,e^{r_1\,s}}n \right) \,\bG\left(\dfrac{y\,e^{r_2\,t_j}}n \right)\,
\\[2mm] \notag
&&\times \PP(\tau_i \in ds\,,\;\sigma_j \in dt_j\,,\; \sigma_k \in dt_k) \leq 2\,\delta\,C^*\,\PP\Bigg( X_i\,e^{-r_1\,\tau_i} > x \,,\; Y_j\,e^{-r_2\,\sigma_j} > y \Bigg)\,,
\eeam
for any $x \wedge y \geq y_0$.

Via symmetrical arguments, we can find similarly that it holds
\beam \label{eq.CKP29V.4.36} 
L_3(x,\,y,\,n) =o\Bigg[ \sum_{i=1}^{n} \sum_{j=1}^{n} \PP\Bigg( X_i\,e^{-r_1\,\tau_i} >x\,,\;Y_j\,e^{-r_2\,\sigma_j} >y\Bigg) \Bigg]\,,
\eeam

Putting \eqref{eq.CKP29V.4.33}, \eqref{eq.CKP29V.4.34} and \eqref{eq.CKP29V.4.36} into  \eqref{eq.CKP29V.4.32} we obtain that for any fixed $n \in \bbn$ it holds
\beam \label{eq.CKP29V.4.37} \notag
\PP\left( \sum_{i=1}^{n} X_i\,e^{-r_1\,\tau_i} > x\,,\;\sum_{j=1}^{n} Y_j\,e^{-r_2\,\sigma_j} >y\right) \lesssim  \sum_{i=1}^{n} \sum_{j=1}^{n} \PP\left( X_i\,e^{-r_1\,\tau_i} > x\,,\; Y_j\,e^{-r_2\,\sigma_j} >y\right)\,.\\
\eeam

From the other hand side, for any fixed $n \in \bbn$, applying Bonferroni inequality twice we find
\beam \label{eq.CKP29V.4.38} \notag 
&&\PP\left( \sum_{i=1}^{n} X_i\,e^{-r_1\,\tau_i} > x\,,\;\sum_{j=1}^{n} Y_j\,e^{-r_2\,\sigma_j} >y\right) \geq \PP\Bigg( \bigcup_{i=1}^{n} \left\{ X_i\,e^{-r_1\,\tau_i} > x \right\}\,, \\[2mm]
&&\;\bigcup_{j=1}^{n} \left\{ Y_j\,e^{-r_2\,\sigma_j} >y \right\}\Bigg)\geq \sum_{i=1}^{n} \PP\left( X_i\,e^{-r_1\,\tau_i} > x\,,\;\bigcup_{j=1}^{n} \left\{ Y_j\,e^{-r_2\,\sigma_j} >y \right\}\right)\\[2mm] \notag
&&- \sum_{1\leq i< k \leq n} \PP\left( X_i\,e^{-r_1\,\tau_i} > x\,,\;X_k\,e^{-r_1\,\tau_k} >x\,,\;\bigcup_{j=1}^{n} \left\{ Y_j\,e^{-r_2\,\sigma_j} >y \right\}\right) \\[2mm] \notag
&&\geq \sum_{i=1}^{n} \sum_{j=1}^{n} \PP\left( X_i\,e^{-r_1\,\tau_i} > x\,,\; Y_j\,e^{-r_2\,\sigma_j} >y \right) \\[2mm] \notag
&&- \sum_{1\leq i< k \leq n} \sum_{j=1}^n\PP\left( X_i\,e^{-r_1\,\tau_i} > x\,,\;X_k\,e^{-r_1\,\tau_k} >x\,,\; Y_j\,e^{-r_2\,\sigma_j} >y \right) \\[2mm] \notag
&&- \sum_{i=1}^{n}  \sum_{1\leq j< k \leq n} \PP\left( X_i\,e^{-r_1\,\tau_i} > x\,,\;Y_k\,e^{-r_2\,\sigma_k} >y\,,\; Y_j\,e^{-r_2\,\sigma_j} >y \right) \\[2mm] \notag
&&\sim \sum_{i=1}^{n} \sum_{j=1}^{n} \PP\left( X_i\,e^{-r_1\,\tau_i} > x\,,\; Y_j\,e^{-r_2\,\sigma_j} >y \right)\,,
\eeam
where we employed similar argument with that for \eqref{eq.CKP29V.4.35}, at the two last sums of the pre-last step. From \eqref{eq.CKP29V.4.37} and \eqref{eq.CKP29V.4.38} we conclude that relation \eqref{eq.CKP29V.4.30} is valid for any fixed $n \in \bbn$.
~\halmos

The following lemma plays crucial role in the proof of Theorem \ref{th.CKP29V.3.2} and extends (under some restrictions) relation \eqref{eq.CKP29V.4.30} uniformly with respect to the summands.

\ble \label{lem.CKP29V.4.5*}
Under the conditions of  Theorem \ref{th.CKP29V.3.2}, asymptotic relation \eqref{eq.CKP29V.4.30} holds uniformly for $n \in \bbn$.
\ele

\pr~ 
At first, from Lemma \ref{lem.CKP29V.4.4} and further, from Lemma \ref{lem.CKP29V.4.3} for any $0 < \delta <1$, there exists some large enough $N \in \bbn$, such that it holds
\beam \label{eq.CKP29V.4.39} \notag 
&&\PP\left( \sum_{i=N+1}^{\infty} X_i e^{-r_1\,\tau_i} > x ,\;\sum_{j=1}^{\infty} Y_j e^{-r_2\,\sigma_j} >y\right)+  \sum_{i=N+1}^{\infty} \sum_{j=1}^{\infty} \PP\left( X_i e^{-r_1\,\tau_i} > x ,\;Y_j e^{-r_2\,\sigma_j} >y \right) \\[2mm]
&&+ \PP\left(  \sum_{i=1}^{\infty}X_i\,e^{-r_1\,\tau_i} > x\,,\;\sum_{j=N+1}^{\infty}Y_j\,e^{-r_2\,\sigma_j} >y \right)\\[2mm] \notag
&&+\sum_{i=1}^{\infty} \sum_{j=N+1}^{\infty} \PP\left( X_i\,e^{-r_1\,\tau_i} > x\,,\;Y_j\,e^{-r_2\,\sigma_j} >y \right)\lesssim \delta\, \PP\left( X_1\,e^{-r_1\,\tau_1} > x\,,\; Y_1\,e^{-r_2\,\sigma_1} >y \right)\,.
\eeam
As \eqref{eq.CKP29V.4.30} holds for any fixed $n \in \bbn$, is implied that it holds uniformly for $1\leq n \leq N$. Hence, it remains to show that \eqref{eq.CKP29V.4.30} holds uniformly for $n > N$. 

For all $n > N$ we obtain
\beam \label{eq.CKP29V.4.40} \notag 
&&\PP\left( \sum_{i=1}^{n} X_i\,e^{-r_1\,\tau_i} > x\,,\;\sum_{j=1}^{n} Y_j\,e^{-r_2\,\sigma_j} >y\right) \geq  \PP\left(  \sum_{i=1}^{N} X_i\,e^{-r_1\,\tau_i} > x\,,\;\sum_{j=1}^{N} Y_j\,e^{-r_2\,\sigma_j} >y \right)\\[2mm] 
&& \sim \sum_{i=1}^{N} \sum_{j=1}^{N} \PP\left( X_i\,e^{-r_1\,\tau_i} > x\,,\;Y_j\,e^{-r_2\,\sigma_j} >y \right) \\[2mm] \notag
&&=\left(\sum_{i=1}^{n} \sum_{j=1}^{n} -\sum_{i=N+1}^{n} \sum_{j=1}^{N} -\sum_{i=1}^{N} \sum_{j=N+1}^{n}  -\sum_{i=N+1}^{n} \sum_{j=N+1}^{n}  \right) \PP\left( X_i\,e^{-r_1\,\tau_i} > x\,,\; Y_j\,e^{-r_2\,\sigma_j} >y \right)\\[2mm] \notag
&&\gtrsim (1-3\,\delta)\sum_{i=1}^{n} \sum_{j=1}^{n} \PP\left( X_i\,e^{-r_1\,\tau_i} > x\,,\; Y_j\,e^{-r_2\,\sigma_j} >y \right)\,,
\eeam 
where at the second step we used \eqref{eq.CKP29V.4.30} for $n=N$, while at the last step we used \eqref{eq.CKP29V.4.39} thrice.

For the compactness of notation we introduce the following
\beao
B:= \left\{ \left(\sum_{i=1}^{N} + \sum_{i=N+1}^{\infty} \right) X_i\,e^{-r_1\,\tau_i} > x\,,\;\left(\sum_{j=1}^{N} + \sum_{j=N+1}^{\infty} \right) Y_j\,e^{-r_2\,\sigma_j} > y \right\}\,.
\eeao

Hence, for all $n > N$ and any $\vep \in (0,\,1)$ it holds
\beam \label{eq.CKP29V.4.41} \notag 
&&\PP\left( \sum_{i=1}^{n} X_i\,e^{-r_1\,\tau_i} > x\,,\;\sum_{j=1}^{n} Y_j\,e^{-r_2\,\sigma_j} >y\right) \leq \PP(B)\\[2mm] \notag
&&=  \PP\left(B\,,\;\left\{  \sum_{i=N+1}^{\infty} X_i\,e^{-r_1\,\tau_i} \leq \vep\,x\right\}\,,\;\left\{\sum_{j=N+1}^{\infty} Y_j\,e^{-r_2\,\sigma_j} \leq \vep\,y\right\} \right)\\[2mm]
&&+  \PP\left(B\,,\;\left\{  \sum_{i=N+1}^{\infty} X_i\,e^{-r_1\,\tau_i} > \vep\,x\right\}\bigcup \left\{\sum_{j=N+1}^{\infty} Y_j\,e^{-r_2\,\sigma_j} >\vep\,y\right\} \right)\\[2mm]  \notag
&& \leq  \PP\left( \sum_{i=1}^{N} X_i\,e^{-r_1\,\tau_i} >(1- \vep)\,x\,,\; \sum_{j=1}^{N} Y_j\,e^{-r_2\,\sigma_j} > (1-\vep)\,y\right)\\[2mm]  \notag
&& + \PP\left(B\,,\; \sum_{i=N+1}^{\infty} X_i\,e^{-r_1\,\tau_i} \leq \vep\,x\,,\; \sum_{j=N+1}^{\infty} Y_j\,e^{-r_2\,\sigma_j} >\vep\,y \right)\\[2mm]  \notag
&& + \PP\left(B\,,\; \sum_{i=N+1}^{\infty} X_i\,e^{-r_1\,\tau_i} > \vep\,x\,,\; \sum_{j=N+1}^{\infty} Y_j\,e^{-r_2\,\sigma_j} \leq \vep\,y \right) \\[2mm]  \notag
&& + \PP\left(B\,,\; \sum_{i=N+1}^{\infty} X_i\,e^{-r_1\,\tau_i} > \vep\,x\,,\; \sum_{j=N+1}^{\infty} Y_j\,e^{-r_2\,\sigma_j} > \vep\,y \right) \\[2mm]  \notag
&& \leq \PP\left( \sum_{i=1}^{N} X_i\,e^{-r_1\,\tau_i} >(1- \vep)\,x\,,\; \sum_{j=1}^{N} Y_j\,e^{-r_2\,\sigma_j} > (1-\vep)\,y\right)\\[2mm]  \notag
&& + \PP\left( \sum_{i=1}^{N} X_i\,e^{-r_1\,\tau_i} >(1- \vep)\,x\,,\; \sum_{j=N+1}^{\infty} Y_j\,e^{-r_2\,\sigma_j} >\vep\,y \right)\\[2mm]  \notag
&& + \PP\left(\sum_{i=N+1}^{\infty} X_i\,e^{-r_1\,\tau_i} > \vep\,x\,,\; \sum_{j=1}^{N} Y_j\,e^{-r_2\,\sigma_j} >(1- \vep)\,y \right) \\[2mm]  \notag
&& + \PP\left( \sum_{i=N+1}^{\infty} X_i\,e^{-r_1\,\tau_i} > \vep\,x\,,\; \sum_{j=N+1}^{\infty} Y_j\,e^{-r_2\,\sigma_j} > \vep\,y \right)=\sum_{i=1}^3 K_i(x,\,y,\,N)\,,
\eeam
For the $K_1(x,\,y,\,N)$ we obtain
\beam \label{eq.CKP29V.4.42} 
&&K_1(x,\,y,\,N) \sim \sum_{i=1}^{N}\,\sum_{j=1}^{N} \,\PP\left(  X_i\,e^{-r_1\,\tau_i} >(1- \vep)\,x\,,\; Y_j\,e^{-r_2\,\sigma_j} > (1-\vep)\,y\right) \lesssim \\[2mm]  \notag
&&\sum_{i=1}^{N}\,\sum_{j=1}^{N} \,\PP\left(  X_i\,e^{-r_1\,\tau_i} >x\,,\; Y_j\,e^{-r_2\,\sigma_j} >  y\right) \leq \sum_{i=1}^{n}\,\sum_{j=1}^{n} \,\PP\left(  X_i\,e^{-r_1\,\tau_i} >x\,,\; Y_j\,e^{-r_2\,\sigma_j} >  y\right)\,,
\eeam
as also $\vep \downarrow 0$, on the second step, where at the first step we used relation \eqref{eq.CKP29V.4.30} for $n=N$, while at the second step we considered that $F,\,G \in \mathcal{C}$, since for any $1 \leq i,\,j \leq N$ via the dominated convergence theorem it holds
\beao
&&\PP\left( X_i\,e^{-r_1\,\tau_i} >(1- \vep)\,x\,,\; Y_j\,e^{-r_2\,\sigma_j} > (1-\vep)\,y\right)\\[2mm] 
&&=\int_0^{\infty} \int_0^{\infty} \bF[(1- \vep)\,x\,e^{r_1\,s}]\,\bG[(1- \vep)\,y\,e^{r_2\,t}] \,\PP\left(  \tau_i \in ds\,,\; \sigma_j \in dt \right)\\[2mm] 
&&=\int_0^{\infty} \int_0^{\infty} \bF(x\,e^{r_1\,s}) \dfrac{\bF[(1- \vep)\,x\,e^{r_1\,s}]}{\bF(x\,e^{r_1\,s})}\,\bG(y\,e^{r_2\,t})\dfrac{\bG[(1- \vep)\,y\,e^{r_2\,t}]}{\bG(y\,e^{r_2\,t})} \,\PP\left(  \tau_i \in ds\,,\; \sigma_j \in dt \right)\\[2mm]  
&&\lesssim \int_0^{\infty} \int_0^{\infty} \bF(x\,e^{r_1\,s})\,\bG(y\,e^{r_2\,t})\,\PP\left(  \tau_i \in ds\,,\; \sigma_j \in dt \right) \\[2mm]  \notag
&&=\PP\left( X_i\,e^{-r_1\,\tau_i} >x\,,\; Y_j\,e^{-r_2\,\sigma_j} > y\right)\,,
\eeao
as also $\vep \downarrow 0$, on the pre-last step. For the $K_2(x,\,y,\,N)$ we find
\beam \label{eq.CKP29V.4.43} \notag 
&&K_2(x,\,y,\,N) \lesssim \delta\,\PP\left( X_1\,e^{-r_1\,\tau_1} >(1- \vep)\,x\,,\; Y_1\,e^{-r_2\,\sigma_1} > \vep\,y\right) \\[2mm] 
&&= \delta\, \int_0^{\infty} \int_0^{\infty}\bF[(1- \vep)\,x\,e^{r_1\,s}]\,\bG(\vep\,y\,e^{r_2\,t}) \,\PP\left(  \tau_1 \in ds\,,\; \sigma_1 \in dt \right)  \\[2mm] \notag
&&\leq \delta\, \int_0^{\infty} \int_0^{\infty} C_2\,\vep^{-q_2}\,\bF(x\,e^{r_1\,s}) \dfrac{\bF[(1- \vep)\,x\,e^{r_1\,s}]}{\bF(x\,e^{r_1\,s})}\,\bG(y\,e^{r_2\,t})\,\PP\left(  \tau_1 \in ds\,,\; \sigma_1 \in dt \right) \\[2mm] \notag
&&\leq \delta\, \int_0^{\infty} \int_0^{\infty} C_2^2\,\vep^{-q_2}\,(1-\vep)^{-p_2}\,\bF(x\,e^{r_1\,s}) \,\bG(y\,e^{r_2\,t})\,\PP\left(  \tau_1 \in ds\,,\; \sigma_1 \in dt \right)  \\[2mm] \notag
&&= \delta\, C_2^2\,\vep^{-q_2}\,(1-\vep)^{-p_2}\,\PP\left(X_1\,e^{-r_1\,\tau_1} >\,x\,,\; Y_1\,e^{-r_2\,\sigma_1} >  y\right) \\[2mm] \notag
&&\leq \delta\, C_2^2\,\vep^{-q_2}\,(1-\vep)^{-p_2}\,\sum_{i=1}^{n}\,\sum_{j=1}^{n}\PP\left(X_i\,e^{-r_1\,\tau_i} > x\,,\; Y_j\,e^{-r_2\,\sigma_j} >  y\right) \,,
\eeam
where at the first step we used relation \eqref{eq.CKP29V.4.39}, while at the third and the fourth steps we used relation \eqref{eq.CKP29V.2.4}. Via symmetric arguments we get
\beam \label{eq.CKP29V.4.44}  
K_3(x,\,y,\,N) \lesssim \delta\, C_2^2\,\vep^{-p_2}\,(1-\vep)^{-q_2}\,\sum_{i=1}^{n}\,\sum_{j=1}^{n}\PP\left(X_i\,e^{-r_1\,\tau_i} > x\,,\; Y_j\,e^{-r_2\,\sigma_j} >  y\right)  \,,
\eeam
Finally, from \eqref{eq.CKP29V.4.39} and \eqref{eq.CKP29V.2.4} we conclude that
\beam \label{eq.CKP29V.4.45}  
&&K_4(x,\,y,\,N) \lesssim \delta\,\PP\left(X_1\,e^{-r_1\,\tau_1} >\vep\,x\,,\; Y_1\,e^{-r_2\,\sigma_1} > \vep\,y \right)   \\[2mm] \notag
&&= \delta\,\int_0^{\infty} \int_0^{\infty} \bF(\vep\,x\,e^{r_1\,s})\,\bG(\vep\,y\,e^{r_2\,t})\,\PP\left(  \tau_1 \in ds\,,\; \sigma_1 \in dt \right)   \\[2mm] \notag
&&= \delta\,\int_0^{\infty} \int_0^{\infty} C_2^2\,\vep^{-(p_2+q_2)}
\bF(x\,e^{r_1\,s})\,\bG(y\,e^{r_2\,t})\,\PP\left(  \tau_1 \in ds\,,\; \sigma_1 \in dt \right)  \\[2mm] \notag
&&= \delta\, C_2^2\,\vep^{-(p_2+q_2)} \PP\left(X_1\,e^{-r_1\,\tau_1} > x\,,\; Y_1\,e^{-r_2\,\sigma_1} >  y\right) \\[2mm] \notag
&&\leq \delta\, C_2^2\,\vep^{-(p_2+q_2)} \,\sum_{i=1}^{n}\,\sum_{j=1}^{n} \PP\left(X_i\,e^{-r_1\,\tau_i} > x\,,\; Y_j\,e^{-r_2\,\sigma_j} >  y\right)\,,
\eeam
Letting in relations \eqref{eq.CKP29V.4.43} - \eqref{eq.CKP29V.4.45} at first $\delta \downarrow 0$ and further $\vep \downarrow 0$ we find
\beam \label{eq.CKP29V.4.46}   \notag
&&K_2(x,\,y,\,N)+K_3(x,\,y,\,N)+K_4(x,\,y,\,N) \\[2mm]
&&=o\left[\sum_{i=1}^{n}\,\sum_{j=1}^{n} \PP\left(X_i\,e^{-r_1\,\tau_i} > x\,,\; Y_j\,e^{-r_2\,\sigma_j} >  y\right)\right]\,,
\eeam
Putting the \eqref{eq.CKP29V.4.42}, \eqref{eq.CKP29V.4.46} into \eqref{eq.CKP29V.4.41} we have that for all $n>N$ it holds
\beam \label{eq.CKP29V.4.47} \notag 
\PP\left(\sum_{i=1}^{n}\, X_i\,e^{-r_1\,\tau_i} > x\,,\; \sum_{j=1}^{n}Y_j\,e^{-r_2\,\sigma_j} >  y\right)\lesssim \sum_{i=1}^{n}\,\sum_{j=1}^{n} \PP\left(X_i\,e^{-r_1\,\tau_i} > x\,,\; Y_j\,e^{-r_2\,\sigma_j} >  y\right)\,,\\
\eeam
which in combination with \eqref{eq.CKP29V.4.40}, letting $\delta \downarrow 0$, gives that \eqref{eq.CKP29V.4.30} holds uniformly for $n>N$.
~\halmos

The next lemma is necessary for the proof of Corollary \ref{cor.CKP29V.3.2}.

\ble \label{lem.CKP29V.4.6}
Under the conditions of Corollary \ref{cor.CKP29V.3.2}, it holds
\beam \label{eq.CKP29V.4.49}  
\PP\left(\sum_{i=1}^{n}\, X_i\,e^{-r_1\,\tau_i} > x\,,\; \sum_{j=1}^{n}Y_j\,e^{-r_2\,\sigma_j} >  y\right) \sim \bF(x)\,\bG(y) \sum_{i=1}^{n}\, \sum_{j=1}^{n} \E[ e^{-\alpha\,r_1\,\tau_i}\,e^{-\beta\,r_2\,\sigma_j}]\,,
\eeam
uniformly for $n \in \bbn$.
\ele

\pr~
For any $i,\,j \in \bbn$ it holds
\beam \label{eq.CKP29V.4.50} \notag 
&&\PP\left(X_i\,e^{-r_1\,\tau_i} > x\,,\;Y_j\,e^{-r_2\,\sigma_j} >  y\right)\\[2mm]
&&=\int_0^{\infty} \int_0^{\infty}\bF(x\,e^{r_1\,s})\,\bG(y\,e^{r_2\,t})\,\PP(\tau_i \in ds,\,\sigma_j \in dt)\\[2mm] \notag
&&\sim \bF(x)\,\bG(y) \int_0^{\infty} \int_0^{\infty} e^{-\alpha\,r_1\,s}\,e^{-\beta\,r_2\,t}\,\PP(\tau_i \in ds,\,\sigma_j \in dt)= \bF(x)\,\bG(y) \E[ e^{-\alpha\,r_1\,\tau_i}\,e^{-\beta\,r_2\,\sigma_j}]\,,
\eeam
where at the second step we took into account that $F \in \mathcal{R}_{-\alpha}$, $G \in \mathcal{R}_{-\beta}$, through the dominated convergence theorem,  the application of which, is due to relations  \eqref{eq.CKP29V.2.3} and  \eqref{eq.CKP29V.2.4}. From  \eqref{eq.CKP29V.4.30} and  \eqref{eq.CKP29V.4.50}, for any $N \in \bbn$, relation \eqref{eq.CKP29V.4.49} holds uniformly for $1\leq n \leq N$.

It remains to show that \eqref{eq.CKP29V.4.49} holds uniformly for $n > N$. At first, for all $n > N$ we obtain that
\beam \label{eq.CKP29V.4.51} \notag 
&&\PP\left(\sum_{i=1}^{n}\, X_i\,e^{-r_1\,\tau_i} > x\,,\; \sum_{j=1}^{n}Y_j\,e^{-r_2\,\sigma_j} >  y\right)\geq \\[2mm] \notag
&&\PP\left(\sum_{i=1}^{N}\, X_i\,e^{-r_1\,\tau_i} > x\,,\; \sum_{j=1}^{N}Y_j\,e^{-r_2\,\sigma_j} >  y\right)\sim \bF(x)\,\bG(y)\,\sum_{i=1}^{N} \sum_{j=1}^{N} \E[ e^{-\alpha\,r_1\,\tau_i}\,e^{-\beta\,r_2\,\sigma_j}]\\[2mm]
&&= \bF(x)\,\bG(y)\,\left( \sum_{i=1}^{n}- \sum_{i=N+1}^{n} \right) \,\left( \sum_{j=1}^{N}- \sum_{j=N+1}^{n} \right) \,\E[ e^{-\alpha\,r_1\,\tau_i}\,e^{-\beta\,r_2\,\sigma_j}] \\[2mm] \notag
&&= \bF(x)\,\bG(y)\,\left( \sum_{i=1}^{n}\,\sum_{j=1}^{n}-\sum_{i=1}^{n}\, \sum_{j=N+1}^{n} - \sum_{i=N+1}^{n}\,\sum_{j=1}^{n} +  \sum_{i=N+1}^{n}\, \sum_{j=N+1}^{n} \right)\,\E[ e^{-\alpha\,r_1\,\tau_i}\,e^{-\beta\,r_2\,\sigma_j}] \\[2mm] \notag
&&= E_{11}(x,\,y;\,N)-E_{12}(x,\,y;\,N)-E_{21}(x,\,y;\,N)+E_{22}(x,\,y;\,N)\,,
\eeam
where at the second step we used relation \eqref{eq.CKP29V.4.49} for $n=N$.

At first we observe that
\beam \label{eq.CKP29V.4.52} 
E_{11}(x,\,y;\,N)=\bF(x)\,\bG(y)\,\sum_{i=1}^{n}\,\sum_{j=1}^{n}\,\E[ e^{-\alpha\,r_1\,\tau_i}\,e^{-\beta\,r_2\,\sigma_j}]\,,
\eeam

For $E_{12}(x,\,y;\,N)$ we find
\beao
&&E_{12}(x,\,y;\,N) \leq \bF(x)\,\bG(y)\,\sum_{i=1}^{n}\,\E^{1/2}[ e^{-2\,\alpha\,r_1\,\tau_i}]\,\sum_{j=N+1}^{\infty}\,\E^{1/2}[e^{-\beta\,r_2\,\sigma_j}]  \\[2mm] \notag
&&\leq \bF(x)\,\bG(y)\,\sum_{i=1}^{\infty}\,\E^{1/2}[ e^{-2\,p_1\,r_1\,\tau_i}] \,\sum_{j=N+1}^{\infty}\,\E^{1/2}[ e^{-2\,q_1\,r_2\,\sigma_j}]\,,
\eeao
where at the first step we used H\"{o}lder inequality. By Assumption  \ref{ass.CKP29V.3.2}, we have that for any $\delta >0$, there exists some large enough $N \in \bbn$, such that it holds
\beam \label{eq.CKP29V.4.53} 
E_{12}(x,\,y;\,N)\leq \delta \,E_{11}(x,\,y;\,N)\,,
\eeam 
and through symmetric argumentation we also obtain
\beam \label{eq.CKP29V.4.54} 
E_{21}(x,\,y;\,N)\leq \delta \,E_{11}(x,\,y;\,N)\,.
\eeam 
Finally, again via  H\"{o}lder inequality and Assumption  \ref{ass.CKP29V.3.2}, we find
\beam \label{eq.CKP29V.4.55}  \notag
&&E_{22}(x,\,y;\,N)\leq 
\bF(x)\,\bG(y)\,\sum_{i=N+1}^{\infty}\,\E^{1/2}[ e^{-2\,p_1\,r_1\,\tau_i}] \,\sum_{j=N+1}^{\infty}\,\E^{1/2}[ e^{-2\,q_1\,r_2\,\sigma_j}] \\[2mm]
&& \leq \delta \,E_{11}(x,\,y;\,N)\,.
\eeam
Putting relations \eqref{eq.CKP29V.4.52} - \eqref{eq.CKP29V.4.55} into \eqref{eq.CKP29V.4.51} and letting $\delta \downarrow 0$,  for all $n > N$ we obtain
\beam \label{eq.CKP29V.4.56}  
\PP\left(\sum_{i=1}^{n}\, X_i\,e^{-r_1\,\tau_i} > x\,,\; \sum_{j=1}^{n}Y_j\,e^{-r_2\,\sigma_j} >  y\right)\gtrsim \bF(x)\,\bG(y)\,\sum_{i=1}^{n}\,\sum_{j=1}^{n}\,\E[ e^{-\alpha\,r_1\,\tau_i}\, e^{-\beta\,r_2\,\sigma_j}] \,.
\eeam

From the other hand side, for all $n>N$ it holds
\beam \label{eq.CKP29V.4.57}  \notag
&&\PP\left(\sum_{i=1}^{n}\, X_i\,e^{-r_1\,\tau_i} > x\,,\; \sum_{j=1}^{n}Y_j\,e^{-r_2\,\sigma_j} >  y\right) \leq \PP\left(\sum_{i=1}^{\infty}\, X_i\,e^{-r_1\,\tau_i} > x\,,\; \sum_{j=1}^{\infty} Y_j\,e^{-r_2\,\sigma_j} >  y\right)\\[2mm]
&&\sim \left(\sum_{i=1}^{N} + \sum_{i=N+1}^{\infty}\right)\, \left(\sum_{j=1}^{N} +\sum_{j=N+1}^{\infty} \right)\,\PP\left( X_i\,e^{-r_1\,\tau_i} > x\,,\; Y_j\,e^{-r_2\,\sigma_j} >  y\right) \\[2mm] \notag
&&=\left(  \sum_{i=1}^{N} \sum_{j=1}^{N}+ \sum_{i=1}^{N} \sum_{j=N+1}^{\infty}+ \sum_{i=N+1}^{\infty} \sum_{j=1}^{N}+ \sum_{i=N+1}^{\infty} \sum_{j=N+1}^{\infty} \right)\,\PP\left( X_i\,e^{-r_1\,\tau_i} > x\,,\; Y_j\,e^{-r_2\,\sigma_j} >  y\right) \\[2mm] \notag
&&=: J_{11}(x,\,y;\,N) + J_{12}(x,\,y;\,N) + J_{21}(x,\,y;\,N) + J_{22}(x,\,y;\,N)\,,
\eeam
where at the second step we applied Lemma \ref{lem.CKP29V.4.5*} for $n=\infty$.

From \eqref{eq.CKP29V.4.50} we find that for all $n>N$, it holds
\beam \label{eq.CKP29V.4.58}  \notag
&&J_{11}(x,\,y;\,N)\sim \bF(x)\,\bG(y)\,\sum_{i=1}^{N} \sum_{j=1}^{N}\,\,\E[ e^{-\alpha\,r_1\,\tau_i}\, e^{-\beta\,r_2\,\sigma_j}] \\[2mm] 
&&\leq  \bF(x)\,\bG(y)\,\sum_{i=1}^{n}\,\sum_{j=1}^{n}\,\E[ e^{-\alpha\,r_1\,\tau_i}\, e^{-\beta\,r_2\,\sigma_j}]\,.
\eeam
For $J_{12}(x,\,y;\,N)$, from \eqref{eq.CKP29V.4.39}, for any $0<\delta <1$ there exists $N \in \bbn$, such that for all $n > N$ we obtain
\beam \label{eq.CKP29V.4.59} 
&&J_{12}(x,\,y;\,N)\lesssim \delta\,\PP\left( X_1\,e^{-r_1\,\tau_1} > x\,,\; Y_1\,e^{-r_2\,\sigma_1} >  y\right) \\[2mm] \notag
&&\sim \delta\, \bF(x)\,\bG(y)\,\E[ e^{-\alpha\,r_1\,\tau_1}\, e^{-\beta\,r_2\,\sigma_1}] \leq \delta\, \bF(x)\,\bG(y)\,\sum_{i=1}^{n}\,\sum_{j=1}^{n}\,\E[ e^{-\alpha\,r_1\,\tau_i}\, e^{-\beta\,r_2\,\sigma_j}]\,,
\eeam
where at the second step we used \eqref{eq.CKP29V.4.50}. Similarly, with  \eqref{eq.CKP29V.4.59}, we can find that for $m=1,\,2$, for all $n > N$, it holds
\beam \label{eq.CKP29V.4.60} 
J_{1m}(x,\,y;\,N)\lesssim \delta\,\bF(x)\,\bG(y)\,\sum_{i=1}^{n}\,\sum_{j=1}^{n}\,\E[ e^{-\alpha\,r_1\,\tau_i}\, e^{-\beta\,r_2\,\sigma_j}]\,.
\eeam  
Putting the \eqref{eq.CKP29V.4.58} - \eqref{eq.CKP29V.4.60} into \eqref{eq.CKP29V.4.57}, and letting $\delta \downarrow 0$, we conclude that for all $n>N$ it holds
\beam \label{eq.CKP29V.4.61} 
\PP\left(\sum_{i=1}^{n}\, X_i\,e^{-r_1\,\tau_i} > x\,,\; \sum_{j=1}^{n}Y_j\,e^{-r_2\,\sigma_j} >  y\right)\lesssim \bF(x)\,\bG(y)\,\sum_{i=1}^{n}\,\sum_{j=1}^{n}\,\E[ e^{-\alpha\,r_1\,\tau_i}\, e^{-\beta\,r_2\,\sigma_j}]\,.
\eeam

From \eqref{eq.CKP29V.4.56} and \eqref{eq.CKP29V.4.61}, we find that \eqref{eq.CKP29V.4.49} holds uniformly for $n>N$.
~\halmos

\bre \label{rem.CKP29V.4.1}
Lemmas \ref{lem.CKP29V.4.5} -\ref{lem.CKP29V.4.6} establish in some sense  the presence of the multivariate non-linear single big jump principle, in finite and infinite randomly weighted sums with weights $e^{-r_1\,\tau_i}$ and $e^{-r_2\,\sigma_j}$ for $i,\,j \in \bbn$. Such kind of problems, with more general random weights, 
were studied recently by \cite{li:2018}, \cite{yang:chen:yuen:2024}, \cite{konstantinides:passalidis:2025c}. Although the most of their conditions are general, the frame is somehow different. In \cite[Th. 2.3(2)]{yang:chen:yuen:2024} the result is close to that of Lemma \ref{lem.CKP29V.4.5*}, however the convergence there is valid as $x\wedge y \to 
\infty$.
\ere

Now, we proceed to the proofs of  Theorem \ref{th.CKP29V.3.2} and  Corollary \ref{cor.CKP29V.3.2}.

\noindent{\bf Proof of Theorem \ref{th.CKP29V.3.2}.}~
Applying Lemma \ref{lem.CKP29V.4.5*} for $n=\infty$ we obtain
\beam \label{eq.CKP29V.4.48} \notag 
&&\PP\left(D_1(\infty) > x\,,\;D_2(\infty)>  y\right)=\PP\left(\sum_{i=1}^{\infty}\, X_i\,e^{-r_1\,\tau_i} > x\,,\; \sum_{j=1}^{\infty}Y_j\,e^{-r_2\,\sigma_j} >  y\right)\\[2mm]
&&\sim \sum_{i=1}^{\infty}\,\sum_{j=1}^{\infty} \PP\left(X_i\,e^{-r_1\,\tau_i} > x\,,\; Y_j\,e^{-r_2\,\sigma_j} >  y\right)\\[2mm] \notag
&&=\sum_{i=1}^{\infty}\,\sum_{j=1}^{\infty}\int_0^{\infty} \int_0^{\infty}\,\bF(x\,e^{r_1\,s})\,\bG(y\,e^{r_2\,t})\,\PP\left(  \tau_i \in ds\,,\; \sigma_j \in dt \right) \\[2mm] \notag
&&=\int_0^{\infty} \int_0^{\infty}\,\bF(x\,e^{r_1\,s})\,\bG(y\,e^{r_2\,t})\,\E\left[N(ds)\,M(dt) \right]\,.~\halmos
\eeam

\noindent{\bf Proof of Corollary \ref{cor.CKP29V.3.2}.}~
Relation  \eqref{eq.CKP29V.3.6} is implied immediately after application of Lemma \ref{lem.CKP29V.4.6} for $n=\infty$.
~\halmos

\end{document}